\documentclass[final,1p,times]{elsarticle}

\usepackage{amsmath,amssymb,amsthm,latexsym,graphicx}
\def\R{\mathbb R} \def\F{\mathcal F} \def\P{\mathbb P} \def\N{\mathbb N} \def\L{\mathcal L}\def\H{\mathcal H}\def\V{\mathcal V}\def\A{\mathcal A}\def\B{\mathcal B}\def\E{\mathbb E}  
\newtheorem{thm}{Theorem}[section]
\newtheorem {definition}{Definition}[section]
\newtheorem{lem}[thm]{Lemma}

\newtheorem{rem}{Remark}[section]
\newenvironment{profof}{\noindent{\bf Proof of }\,}{\,$\Box$\par\smallskip}

 \def\beq{\begin{equation}} \def\enq{\end{equation}}
\def\beqa{\begin{eqnarray}} \def\enqa{\end{eqnarray}}
\numberwithin{equation}{section}

\usepackage[bookmarks,a4paper,colorlinks=true]{hyperref}

\journal{Journal}

\begin{document}

\begin{frontmatter}



\title{A Framework of BSDEs with Stochastic Lipschtz Coefficients through Time Change}
\author{Hun O}
\author{Mun-Chol Kim}
\author[]{Chol-Kyu Pak \corref{cor}}
\cortext[cor]{Corresponding author}
\ead{pck2016217@gmail.com}
\address{Faculty of Mathematics, Kim Il Sung University, Pyongyang, Democratic People's Republic of Korea}
\begin{abstract}
In this paper, we suggest an effective technique based on time-change for dealing with a large class of backward stochastic differential equations (BSDEs for short) defined in general space whose drivers have stochastic Lipschtz coefficients.
By studying the deep properties of random time change combined with stochastic integral and measure theory, we show the relation between the BSDEs with stochastic Lipschtz coefficients and the ones with deterministic Lipschtz coefficients and stopping terminal time, so they are possible to be exchanged with each other from one type to another.  
In other words, the stochastic Lipschtz condition is not essential in the context of BSDEs with random terminal time. 
Next, we derive various results by applying our technique to some types of BSDEs such as Brownian motion BSDE or Markov chain BSDE. 
\end{abstract}

\end{frontmatter}


\section{Introduction}
\label{sec1}
Since their first introduction by Bismut \cite{Bi} in the linear case and the nonlinear extension by Pardoux and Peng \cite{PP1}, Backward stochastic differential equations (BSDEs for short) have been developed rapidly with various types of generalizations in the last decades.
\par BSDEs are closely connected to finance, optimal control and partial differential equation etc.(\cite{ElPQ,P2, PP2, YZ}).
\par Most of BSDEs are concerned with the case of constant time horizon and the uniformly Lipschtz conditions on driver. In many environments, the Lipschtz condition is too restrictive to be assumed, so much effort have been devoted to relax it (\cite{BCa, BLM, LM, M}).
\par In this context, El Karoui and Huang \cite{ElH} studied the BSDEs with stochastic Lipschtz coefficients driven by a general c\`adl\`ag martingale and those were developed under weaker conditions in \cite{CFS}.
For the Brownian motion BSDEs, there are some papers going in this direction(\cite{AIP, BC, BK, JQQ, PA}). Particularly, in \cite{AIP}, Section 3,  the existence of the measure solution was stated by the way of examining the weak convergence of a sequence of measures which were constructed using the martingale representation and the Girsanov change of measure. Also, the reflected backward stochastic differential equations or backward doubly stochastic differential equations (BDSDEs) with stochastic Lipschtz coefficients  were studied in \cite{Hu, ML, MS, O1, O2, W}. 
\par Recently, the inclusive and generalized BSDEs with jumps were studied in the context of stochastic Lipschtz condition in \cite{PPS}. 
\par Although the details are slightly different, the most techniques for the BSDEs with stochastic Lipschtz conditions are similar to the procedure of BSDEs with Lipschtz conditions.
\par That is, the techniques consist of using martingale representation theorem, obtaining a priori estimates and finally using the fixed-point arguments.
\par Other technique was also used in \cite{C5}, where the Lipschtz approximation to the driver was introduced, some estimates were obtained for the convergence of approximation sequence and finally it was shown that the limit of this sequence is a unique solution.
\par In this paper, we approach the problem differently by indirect method.
The technique is based on time change represented by stochastic Lipschtz coefficients. This time change converts the BSDEs with stochastic Lipschtz condition to the ones with uniformly Lipschtz condition and stopping terminal time on another stochastic basis and these two BSDEs are equivalent in some sense. So, if we know the results of BSDEs with random terminal time and uniformly Lipschtz coefficients, then the results are easily extended to the ones with stochastic Lipschtz coefficients through our framework. In other words, the stochastic Lipschtz condition is not a problem in a setting of BSDEs with random terminal time.
\par We briefly mention that the opposite argument also holds, that is, the randomness of terminal time do not play an essential role under the stochastic Lipschtz condition. During our discussion, if the integrator of the driver is a general continuous increasing process, it is converted to the typically well-known one, that is, the Lebesgue measure by time change.
\par Consequently, if we study only the BSDEs to stopping time with standard conditions - the driver satisfies the uniformly Lipschtz continuity, the integrator of the driver is Lebesgue measure, then the research on BSDEs with general conditions - the driver satisfies the stochastic Lipschtz condition, the integrator of the driver is a continuous increasing process is just a corollary of that.
\par And we apply our technique to the detailed BSDEs and get some improved and new results.
\par The prototype of BSDEs is of course Wiener-type BSDE, so we first apply our framework to the BSDE driven by Brownian motion. Here, we deal with the stochastic monotonicity condition more generally. It is clear that the better results in the setting of random terminal time we make use of, the better results in the stochastic Lipschtz setting are obtained. On the other hands, the BSDEs with random terminal time were well-studied sufficiently in many papers.
\par We note that our results include the comparison theorem. In fact, it is a natural question what the behavior of comparison theorem will be like by time change.
\par Here, we emphasize that the comparison theorem as well as wellposedness for BSDEs are easily extended to the stochastic one by our technique. With respect to the previous results in the setting of stochastic Lipschtz, we guarantee the results under weaker conditions on parameters. Moreover we show some new results in the various settings for BSDEs.
\par In this paper, we also apply our framework effectively to the Markov chain BSDE. 
The smart feature is that the discussion on the case of uniformly Lipschtz condition is just inherited to the case of stochastic Lipschtz condition under the same conditions on volumes.
\par In general, for the wellposedness of BSDEs with stochastic Lipschtz condition, the stronger integrability conditions are required than ones with uniformly Lipschtz condition. The main reason is on the discounting property of the terminal time.
This discounting property is contributed to the exponential integrability conditions of volumes  and these conditions are influenced by the Lipschtz coefficients. In fact, discounting property is inherited from the monotonicity of the driver.
In our framework, the original BSDE with stochastic Lipschtz condition can be shown as the BSDE to stopping time which is time-changed in reverse and the time-independent discounting rate of this BSDE with constant Lipschtz coefficients is preserved while time change is processed. This means that the stronger integrability conditions are still required if we use the results of BSDEs with random terminal time obtained by using the monotonicity condition as the key tool. 
 But for the Markov chain BSDEs, the results of undiscounted BSDEs to stopping time without assuming the monotonicity which was researched by Samuel N. Cohen\cite{C1} make our technique more effective.
By passing through the proposed framework, we get a new version of Markov chain BSDEs in the case where the driver has stochastic Lipschtz coefficients for the first time. We also give an example of the real model described as the Markov chain BSDEs with stochastic Lipschtz condition. 
At the end of the paper, we also show some further uses of time change for the BSDEs.
\par The rest of this paper is organized as follows.
In Section 2, we suggest a general map from the BSDEs with stochastic Lipschtz coefficients to the ones with uniformly Lipschtz coefficients by the technique of time change. We discuss this for BSDEs in general space as in \cite{C4}.
The applications to the Wiener-type BSDEs are shown in Section 3.
We give new results on Markov chain BSDEs in Section 4.
In Section 5, we give some concluding remarks.


Let us introduce some useful notations which are used in this paper.
Let $(\Omega,\F,\P)$ be a probability space with a filtration ${\mathbb F}:=\{\F_t\}_{t\ge0}$ satisfying the usual conditions.
We shall assume that $\F=\F_\infty$ and $\F_0$ is trivial.
\begin{itemize}
\item $\|\cdot\|$ denotes the standard Euclidean norm. If $z$ is a matrix, $\|z\|$=Trace$[zz^{\mathsf T}]$, where $[\cdot]^{\mathsf T}$ means the vector transpose.
\item $\B(0,\infty)$ denotes the Borel-$\sigma-$ field given on $(0,\infty)$. 
\item $(\overline\Omega,\overline\F)$ means the product measurable space. That is $\overline{\Omega}:=\Omega\times(0,\infty)$ and $\overline{\F}:=\F\times\B(0,\infty)$.
\item $dQ/d\mu$ denotes the Radon-Nikodym derivative of $Q$ with respect to $\mu$, where $Q$ is absolutely continuous with respect to $\mu$.  
If $\mu$ is Lebesgue meausre and $Q$ is the meausre generated by an absolutely continuous function $f$, then we use $f'$ rather than $dQ/d\mu$.
\item $\E^Q[\cdot]$ means the expectation under measure $Q$.  
\item $L^2(\Omega,\F,\P)$ is the space of square-integrable random variables.
\item $\L$ and $\L^{c}$ are the spaces of local martingales and continuous local martingales, respectively.
\item $\H^2$ is the space of square-integrable martingales.
\item $\H_T^2$ is the space of square-integrable martingales on $[0,T]$.
\item $\H_{loc}^2$ is the space of locally square integrable martingales.
\item $L^2 (M):=\bigl\{ Z\ \big|\ Z$ is predictable$, \E\bigl[\int_0^{\infty}\|Z_t\|^2 d<M>_t\bigr]<+\infty\bigr\}$ where $M\in \H^2$.
\item $L_T^2 (M):=\bigl\{ Z\ |\  Z\cdot I_{[0,T]}\in L^2(M)\bigr\}$, where $M\in \H_T^2$.
\item $L_{loc}^2 (M)$ is the space of predictable processes $Z$ for which there exists a localizing sequence $(\tau_n)$ such that 
\[
\E\biggl(\int_0^{\infty}\|Z\|^2 d<M^{\tau _n}>\biggr)=\E\biggl(\int_0^{\tau _n}\|Z\|^2 d<M>\biggr)<+{\infty},
\] where $M\in \H_{loc}^2$.
\item $L_{t,loc}^2 (M):=\bigl\{X\ |\ X\in L_T^2 (M)$ for any $T<\infty\bigl\}$, where $M\in\H_{loc}^2$.
\item $U_T^2:=\biggl{\{}Y\ \big|\ Y$ is c\`adl\`ag, adapted and $\E\biggl[\sup_{t\in[0,T]}\|Y_t\|^2\biggr]<+\infty\biggr{\}}$.
\item $\V$ is the space of c\`adl\`ag, adapted processes which have finite variation on every finite interval.
\item $\V^+:=\{v\in \V\ |\ v$ is increasing$\}$.
\item $\A:=\bigl{\{}A\in\V\ |\ \E\bigl [\textbf{Var}\bigl ( A(\infty)\bigr ) \bigr ]<\infty)\bigr{\}}$.
\item $\A_{loc}$ is the space of processes locally belonging to $\A$, that is the space of processes $X$ for which there exists a localizing sequence $(\tau_n)$ such that $X^{\tau ^n}\in \A$ for all $n$.
\item $\A_{loc}^+:=\bigl \{ X\in \A_{loc}^+\ |\ X$ is increasing$\bigr \}$.
\item $L_{\theta}^2 (0,\tau;\phi):=\biggl{\{} X\ \big|\ X$ is progressive, $\E\biggl[\int_0^\tau \exp(\theta \phi(s))\|X(s)\|^2 ds\biggr]<\infty \biggr{\}}$, \par where $\theta\in\R$, $\tau$ is stopping time and $\phi$ is an increasing process.
\par If $\phi(t)=t$, we write in $L_{\theta}^2 (0,\tau)$.
\item $L_{\theta}^2 (\tau;\phi)$ is the space of random variables $\xi$ such that $\E[\exp(\theta\phi(\tau))\xi^2]<\infty$.
\item $L_{\theta}^{2,\beta}(0,\tau;\phi):=\{Y\ |\ \beta Y\in L_{\theta}^2 (0,\tau;\phi)\}$.
\item $U_{\theta}^2 (0,\tau;\phi):=\{Y\ |\ Y$ is progressive, $\E[\sup\{\exp(\theta \phi(s))\|Y(s)\|^2 :0\le s\le
\tau\}]<\infty\}$ \par If $\phi(t)=t$, we write in $U_{\theta}^2 (0,\tau)$.
\item If we need to show the Eclidean image space $V$, we use $\L_{\theta}^{2}(0,\tau;\phi,V)$, $\L_{\theta}^{2,\beta}(0,\tau;\phi,V)$ etc.
\item $M_{\theta}^{2,\beta}(0,\tau;\phi;V_1;V_2):=L_{\theta}^2 (\tau;\phi;V_1)\times L_{\theta}^{2,\beta}(0,\tau;\phi;V_2)$, where $V_1,V_2$ are Euclidean spaces.
\end{itemize}


\subsection{Intdroducing BSDEs in general space}
\label{sec1.1}
 As in \cite{C4}, we seem to construct the BSDEs assuming only the usual properties of the filtration and that $L^2(\Omega,\F,\P)$ is a separable Hilbert space.
 Unless otherwise indicated, we should read all equalities(and inequalities) as "up to a measure-zero set" throughout this paper.
\begin{definition}\label{definition:definition11}
For $v\in\V^+$, let us define the measure $\mu_v$on $(\overline{\Omega},\overline{\F})$ as follows.
\beq\label{eq:eq11}
\mu_v(A):=\E\biggl[\int_0^{\infty}I_A(\omega,t)dv\biggr],\quad A\in\overline{\F}
\enq
where the integral is taken pathwise in a Stieltjes sense.
\par This measure $\mu_v$ is called the measure induced (or generated) by $v$.
\end{definition}
\noindent Note that if $v\in\A_{loc}^+$ then $\mu_v$ gives a $\sigma-$ finite measure on $(\overline{\Omega},\overline{\F})$.
\par We give a simple version of the well-known Martingale representation theorem below (see \cite{DV} or \cite{D}).
\begin{thm}[Martingale representation theorem]\label{thm:thm11}
Suppose that $L^2(\Omega,\F,\P)$ is a separable Hilbert space with an inner product $X\cdot Y=\E[XY]$.
\\Then there exists a sequence of $\H^2-$martingales, $M=(M^1 ,M^2 ,...)$ such that $<M^i,M^j>=0$ for $i\neq j$ and every $N\in\H^2$ can be represented as
\beq\label{eq:eq12} 
N_t=N_0+\int_0^t Z_u dM_u=N_0+\sum_{i=1}^{\infty}\int_0^t Z_u^i dM_u^i
\enq
 for some sequence of predictable processes, $Z=(Z^1,Z^2,...)$ satisfying $Z\in L^2 (M)$. 
\par And the predictable quadratic variation processes of these martingales $<M^i>$ satisfy \\
$<M^1>\ \succ\ <M^2>\ \succ...,(\succ$ denotes absolute continuity of induced measures$)$. If $(N^i)$ is another such sequence then $<N^i>\cong<M^i>$, where $\cong$ denotes equivalence of induced measures.
\end{thm}
\begin{rem}\label{rem:rem11}
If the space is generated by Brownian motion, the martingale representation theorem holds on infinite interval (see e.g. \cite{D}, Theorem 6 or references therein). This also implies the martingale representation theorem on every finite interval.
\end{rem}
\noindent For a given $k\in\N$, the general type of BSDE is as follows.
\beq\label{eq:eq13}
Y_t=\xi+\int_t^\tau{g(\omega,s,Y_{s-},Z_s)dv_s} -\sum_{i=1}^{\infty}\int_t^{\tau}Z_s^i dM_s^i,
\enq
\\~
where $\tau$ is an $\mathbb{F}-$stopping time, the terminal value $\xi$ is an $\F_{\tau}-$measurable random variable with values in $\R^k$, the driver $g:\Omega\times (0,\infty)\times \R^k\times \R^{k\times{\infty}}\rightarrow\R^k$ is predictable, $v\in\V$ and the integral of driver is the Lebesgue-Stieltjes integral with respect  to the measures generated by the trajectories of $v$.
\par A solution of the BSDE \eqref{eq:eq13} is a pair of processes $(Y,Z)$ taking values in $\R^k\times\R^{k\times{\infty}}$, where $Y$ is progressive and $Z$ is predictable.
\par In this paper, we shall make the follwing assumption on $v$.\\

\textbf{(A0)} $v$ is a continuous and increasing process.\\

It follows from \textbf{(A0)} that $v$ is locally bounded and $v\in\A_{loc}^+$.
\par Noting that the predictable quadratic variation process $<M>$ identifies an induced measure on $\overline{\F}$ defined by \eqref{eq:eq11}, suppose that the induced measure $\mu_{<M^i>}$ has the following Lebesgue decomposition.
\beq\label{eq:eq14}
\mu_{<M^i>}=\bar{m}^{i,1}+\bar{m}^{i,2}, \quad i\in\N, 
\enq
\\~
where $\bar{m}^{i,1}$ is absolutely continuous with respect to $\mu_v$ and  $\bar{m}^{i,2}$ is orthogonal to $\mu_v$.
\\~
From the generalized Radon-Nikodym Theorem (e.g. see \cite{Med}, Chapter 3, Proposition 3.49), there exist two processes $m_t^{i,1},m_t^{i,2}$ such that $\mu_{m^{i,1}}=\bar{m}^{i,1}$ and $\mu_{m^{i,2}}=\bar{m}^{i,2}$.
\\~
More precisely $m_t^{i,j}=d\pi_t^j/d\P,j=1,2$, where $\pi_t^j (B):=\bar{m}^{i,j}\bigl((0,t]\times B\bigr),B\in\F$. Thus
\beq\label{eq:eq15}
<M^i>_t=m_t^{i,1}+m_t^{i,2}.
\enq
We can consider \eqref{eq:eq15} as the Lebesgue decomposition of $<M^i>$.
\\~
Let us introduce the stochastic semi-norm $\|\cdot\|_{M_t}$ which is defined as 
\beq\label{eq:eq16}
\|z_t\|_{M_t}^2 :=\sum_i\biggl[\|z_t^i\|^2 \cdot (d\bar{m}_t^{i,1}/d\mu _v)\biggr]=\sum_i\biggl[\|z_t^i\|^2\cdot (d \mu_{m^{i,1}}/d\mu_v)(\cdot,t)\biggr],
\enq
for every $z_t=(z_t^1,z_t^2,...)\in\R^{k\times\infty}$.
\\~
Now let us consider the finite time BSDE for $T>0$. We give the following result which is a special case of Theorem 6.1 in \cite{C4}. 
\begin{lem}\label{lem:lem 12}
Let $T>0$, $\xi\in L^2 (\Omega,\F_T,\P;\R^k)$ and suppose that $v$ is a deterministic continuous, increasing function which assigns the positive measure to every non-empty interval in $\R^+$. Let $g:\Omega\times [0,T]\times \R^k \times \R^{k\times {\infty}}\rightarrow\R^k$ be a predictable process such that 
\vskip 0.3cm
1.\quad $\E\biggl[\int_0^T \|g(\omega,t,0,0)\|^2 dv_t\biggr]<\infty$
\vskip 0.3cm
2.\quad For any $y,y'\in \R^k$ and $z,z'\in\R^{k\times \infty}$, there exists $c>0$ such that 
\[
\|g(\omega,t,y,z)-g(\omega,t,y',z')\|^2 \le c\bigl[\|y-y'\|^2+\|z-z'\|_{M_t}^2\bigr],\quad dv\times d\P -a.s.
\]
Then the following BSDE has a unique solution in $U_T^2\times L_T^2 (M)$.
\beq\label{eq:eq17}
Y_t=\xi+\int_t^T g(\omega,s,Y_{s-},Z_s)dv_s -\sum_{i=1}^{\infty}\int_t^T Z_s^i dM_s^i.
\enq
\end{lem}
\noindent In the above lemma, the terminal time is constant. We can also consider the BSDE \eqref{eq:eq13} with stopping terminal time. Perhaps the Lipschtz condition on driver will be still essential and there will be some further conditions related to stopping terminal time for the existence and uniqueness of \eqref{eq:eq13}.
We will not do research of the existence and uniqueness of such BSDEs with random terminal time in this paper.
Our main objective is to show a technique by which the results with respect to stochastic Lipschtz condition are derived from the results with respect to the random terminal time which is considered to be already given. 
\section{Time change and BSDEs}\label{sec2}
We begin with the definition of time change (\cite{RY}, Chapter V).
\begin{definition}
A time change $C$ is a family $\{C(s)\ |\ s>0\}$ of stopping times such that the maps $s\rightarrow C(s)$ are almost surely increasing and right continuous.
\end{definition}
\begin{definition}
If $C$ is a time change, a process $X$ is said to be $C-$continuous if $X$ is constant on each interval $[C_{t-},C_t]$. 
\end{definition}
\noindent We can define the stopped $\sigma-$field $\widetilde{\F}_t:=\F_{C(t)}$ and get the new stochastic basis $(\Omega,\F,\P,\widetilde{\mathbb{F}}=\{\widetilde{\F}_t\}_{t\ge 0})$. It can be easily seen that $\widetilde{\mathbb F}$ also satisfies the usual conditions from the property of stopped $\sigma-$fields. If $X$ is $\F-$progressive then $\widetilde{X}_t:=X_{C_t}$ is $\widetilde{\mathbb{F}}-$adapted and the process $\widetilde{X}_t$ is called the time changed process of $X$. We show a typical example of time change below. 
\par Let us consider an increasing and right-continuous adapted process $A$ (so, progressive) with which we associate 
\beq\label{eq:eq21}
C(s):=\inf\{t\ |\ A(t)>s\},
\enq
\noindent where $\inf(\empty)=+\infty$. This process $C(s)$ is called the inverse of $A(s)$ and we write in $A^{-1}(s)$.
\par As the stohcastic basis satisfies the usual conditions and $A$ is progressive, $A^{-1}(s)$ which is the hitting time of $(s,\infty)$ is a stopping time for every $s>0$. And obviously it is increasing and right continuous. Thus $C=A^{-1}=\{A^{-1}(s)|s>0\}$ is a time change.
\par Throughout this section, we suppose that $C$ is almost surly finite and $C_0=0$ and for any progressive measurable process $X_t$, $\widetilde{X}_t$ means the time changed process of it, unless otherwise indicated. And for the space of processes $V$ with respect to $\mathbb{F}$, $\widetilde{V}$ means the corresponding space with respect to $\widetilde{\mathbb{F}}$. For example, $\widetilde{\L}$ means the space of $\widetilde{\mathbb{F}}-$local martingales.
We give some main results concerning the property of time change under $C-$continuity below.
\begin{lem}\label{lem:lem21} (\cite{RY}, Chapter V, Proposition 1.4).
\\~
Let $C$ be a time change on $(\Omega,\F,\P,\mathbb{F})$. If $h$ is $\mathbb{F}-$progrssive, then $\widetilde{h}$ is $\widetilde{\mathbb F}-$progressive. And if $X$ is a $C-$continuous process of finite variation, then 
\[
\int_0^{C_t}h_udX_u=\int_0^t \widetilde{h}_ud\widetilde{X}_u.
\]
\end{lem}
\begin{lem}\label{lem:lem22} (\cite{RY}, Chapter V, Proposition 1.5)
\\~
If $C$ is a time change on $(\Omega,\F,\P,\mathbb F)$ and $M\in \L^c$ satisfies $C-$continuity, then the following hold.
\vskip 0.3cm
I.\quad $\widetilde{M}\in \widetilde{\L}^c$ and $<\widetilde{M}>=\widetilde{<M>}$
\vskip 0.3cm
II.\quad If $h\in L_{t,loc}^2 (M)$, then $\widetilde{h}\in \widetilde{L}_{t,loc}^2 (\widetilde M)$ and for each $t>0$
\[
\int_0^t \widetilde h _u d\widetilde M _u=\int_0^{C_t}h_u dM_u.
\]
\vskip 0.3cm
Moreover, if $\xi$ is a non-negative random variable, then 
\[
\int_0^{\xi}\widetilde h _ud\widetilde M _u=\int_0^{C_{\xi}}h_u dM_u \quad \P-a.s.\ .
\]
\end{lem}
\noindent
Now we show the property of time change for general locally square-integrable martingales.
\begin{lem}\label{lem:lem23}
If $C$ is a time change on $(\Omega,\F,\P,\mathbb F)$ and $M\in \H_{loc}^2$ is $C-$continuous, then the followings hold.
\vskip 0.3cm
I.\quad $\widetilde{M}\in \widetilde{\H}_{loc}^2$ and $<\widetilde{M}>=\widetilde{<M>}$
\vskip 0.3cm
II.\quad If $h\in L_{t,loc}^2 (M)$, $\widetilde{h}\in \widetilde{L}_{t,loc}^2 (\widetilde M)$ and for each $t>0$
\[
\int_0^t \widetilde h _u d\widetilde M _u=\int_0^{C_t}h_u dM_u
\]
\vskip 0.3cm
Moreover if $\xi$ is a non-negative random variable then 
\[
\int_0^{\xi}\widetilde h _ud\widetilde M _u=\int_0^{C_{\xi}}h_u dM_u \quad \P-a.s.
\]
\end{lem}
\begin{proof}\textbf I. For any $L\in\L$, it is easy to see that $\widetilde{L}\in\widetilde{\L}$ from the optional stopping theorem and $C-$continuity of $M$. 
~\\
As $M\in\H_{loc}^2$, the predictable quadratic variation $<M>$ is in $\A_{loc}^+$ and $M^2-<M>$ is a local martingale from the characterization of $\H_{loc}^2$ martingale (see e.g. \cite{Med}, Chapter 3, Proposition 3.64). Therefore $\widetilde{M^2-<M>}=\widetilde M^2-\widetilde{<M>}$ is an $\widetilde{\mathbb F}-$local martingale. 
~\\
 Let $(\tau_n)$ denote the localizing sequence such that $<M>^{\tau_n}\in\A^+$ for every $n$. 
~\\ 
Then $\widetilde{\tau_n}:=C_{\tau_n}^{-1}=\inf\{t:C_t\ge\tau_n\}$ is an $\widetilde{\mathbb F}-$stopping time for every $n$ and $(\widetilde{\tau_n})$ is a localizing sequence. 
~\\ Noting that $M$ is $C-$continuous if and only if $<M>$ is $C-$continuous (see \cite{RY}, Chapter IV, Proposition 1.13), $<M>$ is constant on $[\tau_n,C_{\widetilde\tau_n}]$. 
~\\
 So $\E\bigl[\widetilde{<M>}(\widetilde\tau _n)\bigr]=\E\bigl[<M>(C_{\widetilde{\tau _n}})\bigr]=\E\bigl[<M>({\tau _n})\bigr]<\infty$. Hence $\widetilde{<M>}\in\widetilde\A_{loc}^+$. 
~\\
And $\widetilde{<M>}$ is also $\widetilde{\mathbb F}-$predictable from the $C-$continuity. Accordingly, using again the characterization of $\H_{loc}^2$ martingale, $\widetilde M\in\widetilde\H_{loc}^2$ and $\widetilde{<M>}=<\widetilde M>$.
\par \textbf II. This is a simple consequence of \textbf I and  Lemma \ref{lem:lem21} together with the relation between stochastic integral and quadratic variation.
\end{proof}
\begin{rem}\label{rem:rem21}
Lemma \ref{lem:lem23} still holds for $\H^2-$martingales under $C-$continuity. That is, if $M$ is $\H^2-$martingale satisfying $C-$continuity, then $\widetilde M\in\widetilde{\H}^2$. In this case we use the characterization of $\H^2-$martingales (e.g. see \cite{Med}, Chapter II, Proposition 2.84) and the same procedure is used for the proof.
\end{rem}

\noindent Now we return to the discussion on BSDE. For the BSDE on which we discuss, the sequence of  $\H^2-$martingales
$M^i\ (i=1,2,...)$ has the martingale representation property on $(\Omega,\F,\P,\mathbb{F})$.
\\~
At this point, the martingale representation on $(\Omega,\F,\P,\widetilde F)$ is naturally expected whereas the time changed processes of $M^i(i=1,2,...)$ are $\widetilde{\H}^2-$martingales under $C-$continuity by  Lemma \ref{lem:lem23}.
\begin{lem}\label{lem:lem24}
Let $C$ be a time change and $\H^2-$martingales $M^i\ (i=1,2,...)$ be $C-$continuous. Then the sequence of $\widetilde{\H}^2-$martingales $(\widetilde M ^i)$ has the martingale representation property for any $\widetilde\H^2-$martingale satisfying  $C^{-1}-$continuity such as in    Theorem \ref{thm:thm11}.
\end{lem}
\begin{proof}Let $\widetilde N$ be an $\widetilde\H^2-$martingale satisfying $C^{-1}-$continuity. Then $\widetilde N_t=\widetilde N_{C^{-1}(C_t)}=N_{C_t}$, where $N_t:=\widetilde N_{C_t^{-1}}$. Obviously, $N\in\H^2$ by   Lemma \ref{lem:lem23}.
Therefore using    Theorem \ref{thm:thm11} and  Lemma \ref{lem:lem23},
\[
\widetilde N_t=N_{C_t}=N_0+\sum_{i=1}^{\infty}\int_0^{C_t}Z_u^i dM_u^i=\widetilde N_0+\sum_{i=1}^\infty \widetilde Z_u^i d\widetilde M _u^i,
\]
for some sequence of $\mathbb F-$predictable processes, $(Z^i)$  satisfying
\[
\E\biggl[\sum_{i=1}^\infty\int_0^\infty (Z_u^i)^2 d<M^i>_u\biggr]<+\infty.
\]
Using Lemma \ref{lem:lem21} and Lemma \ref{lem:lem23} again,
\[
\E\biggl [ \sum_{i=1}^{\infty}\int_0^{\infty}(Z_u^i)^ 2 d<M^i>_u\biggr ]=\E\biggl[\sum_{i=1}^\infty\int_0^\infty (\widetilde Z_u^i)^2 d<\widetilde M ^i>_u\biggr].
\]
\noindent This leads to
\beq\label{eq:eq22}
\E\biggl [ \sum_{i=1}^{\infty}\int_0^{\infty}(\widetilde Z _u^i)^2 d<\widetilde M ^i>_u\biggr ]<+\infty.
\enq
\noindent Hence for any $\widetilde N\in\widetilde H ^2$, there exists a sequence of $\widetilde{\mathbb F}-$predictable processes, $\widetilde Z^i(i=1,2,...)$ satisfying \eqref{eq:eq22} such that
\[
\widetilde N_t=\widetilde N_0+\sum_{i=1}^{\infty}\int_0^t \widetilde Z _u^i d\widetilde M_u^i.
\]
Then by using Lemma \ref{lem:lem23}, we can easily deduce that the martingales $\widetilde M^i ,i=1,2,...$  are mutually orthogonal. The absolute continuity of the induced measures and the uniqueness of the representation are similarly proved.
\end{proof}
\noindent If we know the results for the BSDE \eqref{eq:eq14} with uniformly Lipschtz condition, it is possible to extend to the case where the driver has the stochastic Lipschtz coefficients. This is the main argument in this section.
\par Conveniently, we rewrite the BSDE \eqref{eq:eq14} omitting the index $i$ as follows.
\beq\label{eq:eq23}
Y_t=\xi+\int_t^\tau g(\omega,s,Y_{s-},Z_s)dv_s- \int_t^\tau Z_s dM_s,\quad 0\le t\le \tau,
\enq
\noindent where $Z=(Z^1,Z^2,...)$ and $M=(M^1,M^2,...)$.
\par Assume that the driver of \eqref{eq:eq23} satisfies the following stochastic Lipschtz condition.\\

\textbf{(A1)} There exist predictable processes $r_t$ and $u_t$ such that
\[
\|g(\omega,t,y_t ,z_t)-g(\omega,t,y'_t ,z'_t)\|\le r_t\|y_t -y'_t\|+u_t\|z_t -z'_t\|_{M_t}, \quad d\mu_v-a.s.
\]
 for any $y_t,y'_t\in\R^k$ and $z_t,z'_t\in \R^{k\times \infty}$, where $\alpha_t^2:=\max\{r_t,u_t^2\}>\epsilon$ for some $\epsilon >0$ and $\alpha_t^2$ is pathwise Stieltjes-integrable with respect to $v$ for every finite interval in $\R^+$.\\

 Now we define the following process.
\beq\label{eq:eq24}
\phi (t):=\int_0^t \alpha_s^2 dv_s.
\enq
The remarkable point is that $\phi^{-1}$ i.e. the inverse of $\phi_t$ defined by \eqref{eq:eq21} is a time change. We shall make a good use of this process in the view of time change. It is clear that $\phi^{-1}$ is a.s. finite and $\phi^{-1}(0)=\phi(0)=0$. From now, the symbol $C$ which has meant time change will be replaced by $\phi^{-1}$. 
The focus of this section is on the technique, so we do not have detailed discussion on the space of solutions.
The main result in this section is as follows.
\begin{thm}\label{thm:thm25}
Let $\phi(t)$ be a process defined by \eqref{eq:eq24} and $M$ be $\phi^{-1}-$continuous. If $(Y_t,Z_t)$ is a solution of BSDE \eqref{eq:eq23} satisfying $\textbf{(A0)}$ and $\textbf{(A1)}$ on $(\Omega,\F,\P,\mathbb F)$, then $(y_t,z_t):=(Y_{\phi_t^{-1}(t)},Z_{\phi_t^{-1}(t)})$ is a solution of the following BSDE on $(\Omega,\F,\P,\widetilde{\mathbb F})$.
\beq\label{eq:eq25}
y_t=\xi+\int_t^{\widetilde{\tau}} \widetilde g(\omega,s,y_{s-},z_s)ds-\int_t^{\widetilde{\tau}}z_s d\widetilde M_s,\quad 0\le t\le \widetilde{\tau},
\enq
where
\beq\label{eq:eq26}
\widetilde g(\omega,s,y,z):=g(\omega,\phi^{-1}(s),y,z)/ \alpha^2(\phi^{-1}(s)),\quad \widetilde{\tau}:=\phi(\tau),\quad \widetilde M_s :=M_{\phi_s^{-1}}.
\enq
The converse is also true, that is if $(y_t,z_t)$ is a solution of the BSDE \eqref{eq:eq25}, then $(Y_t,Z_t):=(y_{\phi (t)},z_{\phi (t)})$ is a solution of the BSDE \eqref{eq:eq23}. Mainly the new driver $\widetilde g$ of \eqref{eq:eq25} satisfies uniform Lipschtz continuity such that for any $y_t,y'_t\in\R^k$ and $z_t,z'_t\in \R^{k\times \infty}$,
\[
\|\widetilde g(\omega,t,y_t ,z_t)-\widetilde g(\omega,t,y'_t ,z'_t)\|\le \|y_t -y'_t\|+\|z_t -z'_t\|_{M_t}, \quad dt\times d\P-a.s.\ .
\]
\end{thm}
\begin{proof}
We split the proof into four steps.
\vskip 0.2cm
\textbf{Step 1}
We first show that $\widetilde v(\cdot,\omega):=v\bigl(\phi^{-1}(\cdot,\omega),\omega\bigr)$ is absolutely continuous for each $\omega\in\Omega$. As $v$ is increasing and continuous, $v^{-1}$  defined by \eqref{eq:eq21} is a time change and $v$ is $v^{-1}-$continuous. Therefore by Lemma \ref{lem:lem21}, we can see that
\begin{align}
\phi(t)=\int_0^t \alpha^2 (s)dv_s&=\int_0^{v^{-1}(v_t)}\alpha^2(s)dv_s-\int_t^{v^{-1}(v_t)}\alpha^2(s)dv_s\nonumber\\
&=\int_0^{v^{-1}(v_t)}\alpha^2(s)dv_s=\int_0^{v_t}\alpha^2(v_s^{-1})ds=\biggl[\int_0^\cdot \alpha^2(v_s^{-1})ds\circ v\biggr](t)\nonumber.
\end{align}
Thus, $\widetilde v_t=[v\circ{\phi^{-1}](t)}=\biggl[v\circ v^{-1}\circ\biggl(\int_0^{\cdot}\alpha^2(v_s^{-1})ds\biggr)^{-1}\biggr](t)=\biggl(\int_0^{\cdot}\alpha^2(v_s^{-1})ds\biggr)^{-1}(t)$.
\\
 Noting that $\alpha^2(s)>\epsilon$, $\widetilde v^{-1}(\cdot)=\int_0^{\cdot}\alpha^2(v_s^{-1})ds$ is strictly increasing and absolutely continuous for each $\omega\in\Omega$ and so is the reversed process.
Hence $\widetilde v_t$ (resp. $\mu_{\widetilde v}$) is absolutely continuous with respesct to Lebesgue measure (resp $dt\times d\P$) and 
\[
d\mu_{\widetilde v}/(dt\times d\P)=d\widetilde v_t/{dt}=1/{[\alpha^2\circ v^{-1}\circ v\circ\phi^{-1}](t)}=\alpha^{-2}(\phi_t^{-1}).
\]
In fact, we can see that $\widetilde v_t$ (resp. $\mu_{\widetilde v}$) is equivalent to Lebesgue measure (resp. $dt\times d\P)$. We also mention that $v$ is $\phi^{-1}-$continuous.
\vskip 0.2cm
\textbf{Step 2}
We derive the Lebesgue decomposition of the measure induced by ${<\widetilde M>}$.
First, we show that $m^1$ is a.s. $\phi^{-1}-$continuous.
Suppose that $v$ is a constant on $[a,b]\ (0\le a<b)$. Then for any $c\in[a,b]$ and $B\in\F$, $\bar m([c,b]\times B)=\E[\int_c^b I_B(\omega)\cdot (d\bar m^1/d\mu_v)(\omega,t)dv_t]=0$. Noting that $\bar m([0,t]\times B)=\int_B m_t^1d\P$, 
\[
0=\bar m([c,b]\times B)=\bar m([0,b]\times B)-\bar m([0,c]\times B)=\int_B (m_b^1-m_c^1)d\P.
\]  
Hence $m^1$ is a.s. constant on $[a,b]$. Because $v$ is $\phi^{-1}-$continuous from \textbf{Step 1}, we can see that $m^1$ is a.s.  $\phi^{-1}-$continuous.
Recalling \eqref{eq:eq15} and using  Lemma \ref{lem:lem23}, we obtain(omitting the index $i$)
\beq\label{eq:eq27}
<\widetilde M>_t=\widetilde{<M>}_t=\widetilde m_t^1+\widetilde m_t^2.
\enq
And the continuity of $\phi$ which comes from the continuity of $v$ implies $\phi(\phi_t^{-1})=t$.  
Now we can use  Lemma \ref{lem:lem21} to show 
\begin{align}
\mu_{\widetilde m^1}(A)&=\E\biggl[\int_0^{\infty}I_A (\omega,t)d\widetilde m_t^1\biggr]=
\E\biggl[\int_0^{\infty}I_A(\omega,\phi(\omega,t))dm_t^1\biggr]\nonumber\\
&=\E\biggl[\int_0^{\infty}I_A(\omega,\phi(\omega,t))(d\overline m^1/d\mu_v)(\omega,t)\cdot dv_t\biggr]\nonumber\\
&=\E\biggl[\int_0^{\infty}I_A(\omega,t)[d\overline m^1/d\mu_v](\phi_t^{-1})d\widetilde v_t\biggr]\nonumber
\end{align}
for any $A\in\overline\F$, where $[d\overline m^1/d\mu_v](\phi_t^{-1}):=[d\overline m^1/d\mu_v](\omega,\phi^{-1}(\omega,t))$.\\
Thus $\mu_{\widetilde m^1}\prec\mu_{\widetilde v}$ and $d\mu_{\widetilde m^1}/d\mu_{\widetilde v}=[d\overline m^1/d\mu_v](\phi_t^{-1})$.\\
Noting that $\mu_{\widetilde v}\prec dt\times d\P$ by \textbf{Step 1}, we can deduce $\mu_{\widetilde m^1}\prec dt\times d\P$ and 
\beq\label{eq:eq28}
d\mu_{\widetilde m^1}/(dt\times d\P)=[d\overline m^1/d\mu_v](\phi_t^{-1})\cdot d\widetilde v_t/dt.
\enq
Similarly, $\mu_{m^2}$ is orthogonal to $dt\times d\P$. This shows that \eqref{eq:eq27} is the Lebesgue decomposition of $<\widetilde M>$ with respect to $dt\times d\P$.
\vskip 0.2cm
\textbf{Step 3}. As $(Y_t,Z_t)$ is the solution of \eqref{eq:eq23},
\[
y_t:=Y_{\phi^{-1}(t)}=\xi+\int_{\phi^{-1}(t)}^{\tau}g(\omega,s,Y_{s-},Z_s)dv_s-\int_{\phi^{-1}(t)}^{\tau}Z_s dM_s,\quad 0\le t\le \widetilde\tau.
\]
By  Lemma \ref{lem:lem21} and \textbf{Step 1},
\begin{align}
\int_{\phi^{-1}(t)}^{\tau}g(\omega,s,Y_{s-},Z_s)dv_s &=\int_t^{\widetilde\tau}g(\omega,\phi^{-1}(s),Y_{\phi^{-1}(s)-},Z_{\phi^{-1}(s)})d\widetilde v_s/ds\cdot ds\nonumber\\
&=\int_t^{\widetilde\tau}g(\omega,\phi_s^{-1},Y_{\phi^{-1}(s)-},Z_{\phi^{-1}(s)})\alpha^{-2}(\phi^{-1}(s)) ds\nonumber\\
&=\int_t^{\widetilde\tau}\widetilde g(\omega,s,y_{s-},z_s)ds\nonumber.
\end{align}
By  Lemma \ref{lem:lem23} and $\phi^{-1}-$continuity of $M$,
\[
\int_{\phi^{-1}(t)}^{\tau}Z_s dM_s=\int_t^{\phi(\tau)} Z_{\phi^{-1}(s)} dM_{\phi^{-1}(s)}=\int_t^{\widetilde\tau}z_s d\widetilde M_s.
\]
So we have
\[
y_t=\xi+\int_t^{\widetilde\tau}\widetilde g(\omega,s,y_{s-},z_s)-\int_t^{\widetilde\tau}z_s d\widetilde M_s,\quad 0\le t\le\widetilde\tau.
\]
As $Y_t$ is $\mathbb F-$progressive, $y_t$ is $\widetilde{\mathbb F}-$progressive. Due to the fact that all stochastic integrals are indistinguishable from the stochastic integrals of predictable processes, we can consider $z_t$ is predictable. Accordingly, $(y_t,z_t)$ is a solution of BSDE \eqref{eq:eq25} on $(\Omega,\F,\P,\widetilde{\mathbb F})$. Passing back through the above procedure, the converse argument is trivial.
\vskip 0.2cm
 \textbf{Step 4.}
Finally, we show that $\widetilde g$ satisfies uniform Lipschtz continuity.
It follows from the results in \textbf{Step 2} that 
\[
\begin{aligned}
\|z_t\|_{\widetilde M_t}^2=\|z_t\|^2 \cdot [d\mu_{\widetilde m^1}/(dt\times d\P)]&=\|z_t\|^2\cdot  [d\overline m^1/d\mu_v](\phi_t^{-1})d\widetilde v_t/dt\nonumber\\
&=\alpha^{-2}(\phi_t^{-1})\|z_t\|_{M_u}^2\big|_{u=\phi^{-1}(t)}.\nonumber
\end{aligned}
\]
From the stochastic Lipschtz condition on $g$,
\begin{align}
\|\widetilde g(\omega,t,y_t,z_t)&-\widetilde g(\omega,t,y'_t,z'_t)\|=\|\widetilde g(\omega,\phi^{-1}(t),y_t,z_t)-\widetilde g(\omega,\phi^{-1}(t),y'_t,z'_t)\|\alpha^{-2}(\phi_t^{-1})\nonumber\\
&\le\alpha^{-2}(\phi^{-1}(t))\bigl[r_{\phi^{-1}(t)}\|y_t-y'_t\|+u_{\phi^{-1}(t)}\bigl(\|z_t-z'_t\|_{M_s}\big|_{s=\phi^{-1}(t)}\bigr)\bigr]\nonumber\\
&=\alpha^{-2}(\phi^{-1}(t))\bigl[r_{\phi^{-1}(t)}\|y_t-y'_t\|+u_{\phi^{-1}(t)}\alpha(\phi^{-1}(t))\bigl(\|z_t-z'_t\|_{\widetilde M_t}\bigr]\nonumber\\
&=\frac{r_{\phi^{-1}(t)}}{\max\{r_{\phi^{-1}(t)},u_{\phi^{-1}(t)}^2\}}\|y_t-y'_t\|+\frac{u_{\phi^{-1}(t)}}{\sqrt{\max\{r_{\phi^{-1}(t)},u_{\phi^{-1}(t)}^2\}}}\|z_t-z'_t\|_{\widetilde M_t}\nonumber\\
&\le\|y_t-y'_t\|+\|z_t-z'_t\|_{\widetilde M_t},\quad\quad d\mu_{\widetilde v}-a.s.\nonumber
\end{align}
for any $y_t,y'_t\in\R^k$ and $z_t,z'_t\in \R^{k\times \infty}$.
From \textbf{Step 1}, we know that $\mu_{\widetilde v}$ is equivalent to $dt\times d\P$. So Lipschtz property on $\widetilde g$ holds $dt\times d\P-$a.s.
\end{proof}
\begin{rem}\label{rem:rem22}
If the trajectories of $v$ are strictly increasing,  then $\phi^{-1}$ is strictly increasing and continuous (that is $\phi^{-1}(\phi(t))=\phi(\phi^{-1}(t))=t$), so we do not have to assume that $M$ is $\phi^{-1}-$continuous. Remark that $\phi^{-1}-$ continuity of $M$ is equivalent to $\phi^{-1}-$ continuity of $m^2$. 
\end{rem}
\begin{rem}\label{rem:rem23}
In our discussion, the continuity of $v$ which leads to the continuity of $\phi$, plays an important role. This guarantees $v(v^{-1}(t))=\phi(\phi^{-1}(t))=t$. If $v$ is a finite variation process possibly with jumps, it may be needed to decompose the Stieltjes measures generated by the trajectories of $v$ as the continuous part and the discontinuous one. Perhaps it may be non-trivial.
\end{rem}
\begin{rem}\label{rem:rem24}
If we only want to simplify the continuous integrator of driver, it is sufficient to use $v^{-1}$ as the time change. 
\end{rem}
It is natural to try the comparison theorem under the stochastic Lipschtz condition by means of time change. Suppose that we have two BSDEs satisfying \textbf{(A0)}, \textbf{(A1)} and let $(g,\xi)$,$(\bar g,\bar\xi)$ be the corresponding generators. And let $(Y,Z)$,$(\bar Y,\bar Z)$ be the associated solutions. The following assumption plays an important role to ensure that the comparison theorem holds (\cite{C4}).\\

\noindent \textbf{(A2)} 
\begin{enumerate}
\item  For every $j$, there exists $\hat\P_j$ equivalent to $\P$ such that $j^{th}$ component of $X$ as defined by 
\[
e_j^{\mathsf T}:=-\int_0^r e_j^{\mathsf T}[g(\omega,u,\bar Y_{u-},Z_u)- g(\omega,u,\bar Y_{u-},\bar Z_u)]dv_u+\int_0^r e_j^{\mathsf T}[Z_u-\bar Z_u]dM_u
\]
is $\hat\P_j-$supermartingale.
\item If for all $r\ge 0$, 
\[
e_j^{\mathsf T}Y_r-\E^{\hat\P_i}\biggl[\int_r^\infty e_i^{\mathsf T}g(\omega,u,Y_{u-},Z_u)dv_u\big |\F_r\biggr]\ge e_j^{\mathsf T}\bar Y_r-\E^{\hat\P_i}\biggl[\int_r^\infty e_i^{\mathsf T}g(\omega,u,\bar Y_{u-},Z_u)dv_u\big |\F_r\biggr]\nonumber
\]
for all $i$, then $Y_r\ge \bar Y_r$ for all $r\ge 0$ componentwise.
\end{enumerate}
The driver satisfying \textbf{(A2)} is often called balanced. This notation originated from finance, as in some sense, the driver balances the outcomes to hedge. This driver is closely connected to no-arbitrage opportunity and furthermore the condition under which the comparison theorem holds for martingale-type BSDEs possibly with jumps (see \cite{C3, C4} or \cite{C5}, Part IV).
\par It is obvious that the comparison theorem holds for BSDE \eqref{eq:eq23} if and only if the comparison theorem holds for the corresponding BSDE \eqref{eq:eq25}. Now we shall show that the essential conditions which ensure that the comparison theorem holds are preserved while the time change is processed. We still assume that the BSDE satisfies \textbf{(A0)},\textbf{(A1)} and $M$ is $\phi^{-1}-$continuous. 
\begin{thm}\label{thm:thm26}
If BSDE \eqref{eq:eq23} satisfies \textbf{(A2)}, the time changed BSDE \eqref{eq:eq25} also satisfies \textbf{(A2)} with respect to filtration $\widetilde{\mathbb F}$.
\end{thm}
\begin{proof}
First by the optional stopping theorem, $e_j^{\mathsf T}\widetilde X$ is $\widetilde{\mathbb F}-$supermartingale under $\hat\P_j$ for every $j$ using that $e_j^{\mathsf T}$ is $\mathbb F-$supermartingale under $\hat\P_j$.
Having the similar procedure to \textbf{Step 3} in the proof of    Theorem \ref{thm:thm25}, we obtain  
\[
e_j^{\mathsf T}\widetilde X_r=-\int_0^re_j^{\mathbb T}[\widetilde g(\omega,u,\bar y_{u-},z_u)-\widetilde g(\omega,u,\bar y_{u-},\bar z_u)]du+\int_0^r e_j^{\mathbb T}[z_u-\bar z_u]d\widetilde M_u
\]
So the first part of \textbf{(A2)} is satisfied with respect to $\widetilde{\mathbb F}$ for BSDE \eqref{eq:eq25}. Similarly we can prove that the second part is also satisfied. 
\end{proof} 
We conclude this section with the following statement.
\vskip 0.2cm
\textbf{Interesting remark on terminal time}
~\\
When we study the BSDEs with stochastic Lipschtz coefficients, the randomness of terminal time does not play an important role.
This is illustrated as follows.
Due to the Remark \ref{rem:rem24}, we can suppose that the BSDE is given in the following type without loss of generality.
\beq\label{eq:eq29}
Y_t=\xi+\int_t^{\tau}g(\omega,s,Y_{s-},Z_s)ds-\int_t^{\tau}Z_s dM_s.
\enq
We use the following process introduced for the quadratic BSDEs in \cite{AIP}.
\[
\Phi(\omega,t):=\frac{t}{1+\tau\wedge t},\quad t\ge 0.
\]
After the simple calculation, we get
\[
\Phi^{-1}(t)=t/(1-t),\quad [\Phi^{-1}(t)]'=-(t-1)^{-2},\quad 0\le t\le \widetilde{\tau}=\Phi(\tau)<1.
\]
Obviously $\Phi^{-1}$ is time change and we can deduce the following BSDE on $(\Omega,\F,\P,\mathbb F)$ equivalent to \eqref{eq:eq29} in some sense.
\beq\label{eq:eq210}
y_t=\xi+\int_t^1 G(\omega,s,y_{s-},z_s)ds-\int_t^1 z_sd\widetilde M_s,
\enq
where $G(\omega,s,y,z):=I_{s\le \tau}g(\omega,\Phi^{-1}(s),y,z)[\Phi^{-1}(s)]'$ and $\widetilde M_s:=M_{\Phi^{-1}(s)}$.
We mention that the new driver $G$ is stochastic Lipschtz even though the original driver $g$ is uniform Lipschtz.
In fact, if we suppose that $g$ has constants $r,u$ as the Lipschtz coefficients, for any $y,y'\in\R^k$ and $z,z'\in\R^{k\times\infty}$,
\begin{align}
\|G(\omega,s,y,z)-G(\omega,s,y',&z')\|=I_{s\le\widetilde{\tau}}|(\Phi^{-1})'(s)|\cdot\|g(\omega,\Phi_s^{-1},y,z)-g(\omega,\Phi_s^{-1},y',z')\|\nonumber\\
&\le I_{s\le\widetilde{\tau}}|(\Phi^{-1})'(s)|[r\|y-y'\|+u\|z-z'\|_{\widetilde M_s}|(\Phi^{-1})'(s)|^{-1/2}]\nonumber\\
&=I_{s\le\widetilde{\tau}}[r(1-s)^{-2}\|y-y'\|+u(1-s)^{-1}\|z_t-z'_t\|_{\widetilde M_s}]\nonumber\\
&\le r(1-\widetilde{\tau})^{-2}\|y-y'\|+u(1-\widetilde\tau)^{-1}\|z-z'\|_{\widetilde M_s}\nonumber\\
&=r(1+\tau)^2\|y-y'\|+u(1+\tau)\|z-z'\|_{\widetilde M_s}\nonumber.
\end{align}
This means that the stopping terminal time of BSDEs can be converted to constant and this operation is adapted to the class of BSDEs with stochastic Lipschtz condition.
\section{Wiener-type BSDEs with stochastic monotone coefficients}\label{sec3}
The well-known and mostly studied type of BSDEs are of course Wiener-type BSDEs.
Let $W$ be $d$-dimensional Brownian motion on $(\Omega,\F,\P)$ and $\mathbb F:=\{\F_t\}_{t\ge 0}$ be the natural complete, right continuous filtration generated by $W$.
It is worthy to study Wiener-type BSDEs with stochastic Lipschtz conditions. 
For example, let us consider the pricing problem of a European contingent claim. This problem is equivalent to solving the following linear BSDE:
\[
Y_t=\xi+\int_t^T(r_s Y_s+u_s Z_s)ds+\int_t^T Z_s dW_s,
\]
where $\xi$ is the contingent claim to hedge, $r_s$ is the interest rate, $u_s$ is the risk premium vector and $T$ is the maturity date.
In general, $r_t$ and $u_t$ both will not be bounded and moreover the maturity date will be non-deterministic. In this case the Lipschtz condition does not hold uniformly any more. For the forward-backward BSDEs, when the uncertainty of driver only comes from a solution of forward component, we can give the probabilistic interpretation of a system of semi-linear elliptic PDEs (see \cite{BEN}, Remark 4.6). 
\par We shall have slightly different procedure from Section \ref{sec2}, but this is essentially the same.
For the discussion of Martingale-type BSDE, the martingale term is changed into a martingale on another stochastic basis. 
As the quadratic variations of them are different, these martingales are not equal in general. In view of general Martingale-type BSDE, this is non-sense.
But the theory of Wiener-type BSDE is well studied than others so it will be convenient for the research if the Wiener-type BSDE is converted to the Wiener-type BSDE on a new basis.
We consider the following BSDE driven by Brownian motion on stochastic basis $(\Omega,\F,\P,\mathbb F)$.
\beq\label{eq:eq31}
Y_t=\xi+\int_t^\tau f(\omega,s, Y_s, Z_s)ds+\int_t^\tau Z_s dW_s,\quad 0\le t \le \tau,
\enq
where $\tau$ is an a.s. finite $\mathbb F-$stopping time, $\xi$ is an $\F_{\tau}-$measurable random variable with values in $\R^k$ and $f:\Omega\times\R^+\times\R^k\times\R^{k\times d}\rightarrow\R^k$ is $\mathbb F-$progressive. Due to $<W_t>=t$, the stochastic semi-norm defined by \eqref{eq:eq16} is obtained as $\|z\|_{W_t}=\|z\|$ for $z\in\R^{k\times d}$. Let the driver $f$ satisfy the stochastic Lipschtz condition. That is there exist non-negative progressive processes $r_t$ and $u_t$ such that 
\beq\label{eq:eq32}
\|f(\omega,t,y,z)-f(\omega,t,y',z')\|\le r_t\|y-y'\|+u_t\|z-z'\|,\quad dt\times d\P-a.s.
\enq
for any $y,y'\in\R^k$,\ $z,z'\in\R^{k\times d}$. As in  Section \ref{sec2}, we assume that there exists $\epsilon>0$ such that $\alpha^2(t):=\max\{r_t,u_t^2\}>\epsilon$ and that $\alpha^2(t)$ is Lebesgue-integrable on every finite interval in $\R^+$. And we introduce the following strictly increasing, absolutely continuous process:
\beq\label{eq:eq33}
\phi(t):=\int_0^t\alpha^2(s)ds.
\enq
We set $\widetilde{\F}:=\F_{\phi^{-1}(t)}$. Now we define stochastic process $\widetilde W_t$ as follows.
\beq\label{eq:eq34}
\widetilde W_t:=\int_0^t [(\phi^{-1})'(s)]^{-1/2}dW_{\phi^{-1}(s)}.
\enq
Then $\widetilde W$ is a continuous $\widetilde {\mathbb F}-$local martingale and for each $i$, $<\widetilde W^i>_t=\int_0^t{\bigl[(\phi^{-1})'(s)\bigr]}^{-1}d<W^i>_{\phi^{-1}(s)}=\int_0^t{\bigl[(\phi^{-1})'(s)\bigr]}^{-1}d\phi^{-1}(s)=t$, so it is $\widetilde{\mathbb F}-$Brownian motion by L\`evy's characterization theorem.
 If $(Y_t,Z_t)$ is a solution of \eqref{eq:eq31}, then
\begin{align}
y_t:&=Y_{\phi^{-1}(t)}=\xi+\int_{\phi^{-1}(t)}^{\tau}f(\omega,s,Y_s,Z_s)ds-\int_{\phi^{-1}(s)}^{\tau}Z_s dW_s\nonumber\\
&=\xi+\int_t^{\widetilde{\tau}}f(\omega,\phi^{-1}(s),Y_{\phi^{-1}(s)},Z_{\phi^{-1}(s)})(\phi^{-1})'(s)ds-\int_t^{\widetilde{\tau}}Z_{\phi^{-1}(s)} dW_{\phi^{-1}(s)}\nonumber\\
&=\xi+\int_t^{\widetilde{\tau}}f(\omega,\phi^{-1}(s),y_s,Z_{\phi^{-1}(s)})(\phi^{-1})'(s)ds-\int_t^{\widetilde{\tau}}Z_{\phi^{-1}(s)} {\big[(\phi^{-1})'(s)\bigr]}^{1/2} d\widetilde W_s\nonumber\\
&=\xi+\int_t^{\widetilde{\tau}}f\bigl(\omega,\phi^{-1}(s),y_s,z_s {\bigl[(\phi^{-1})'(s)\bigr]}^{-1/2}\bigr)(\phi^{-1})'(s)ds-\int_ t^{\widetilde{\tau}} z_s d\widetilde W_s,\ 0\le t\le\widetilde{\tau}\nonumber,
\end{align}
where $z_s:=Z_{\phi^{-1}(s)}{\bigl[(\phi^{-1})'(s)\bigr]}^{1/2}$ and  $\widetilde\tau:=\phi(\tau)$.
So, if we set $\widetilde f$ as
\beq\label{eq:eq35}
\widetilde f(\omega,s,y,z):=f\bigl(\phi^{-1}(s),y,z\cdot \bigl[(\phi^{-1})'(s)\bigr]^{-1/2}\bigr)\cdot(\phi^{-1})'(s)
\enq
then $(y_t,z_t)=(Y_{\phi^{-1}(t)},Z_{\phi^{-1}(t)}\cdot (\phi^{-1})'(t)^{1/2})$ is the solution of the following BSDE on $(\Omega,\F,\P,\widetilde{\mathbb F})$:
\beq\label{eq:eq36}
y_t=\xi+\int_t^{\widetilde\tau}\widetilde f(\omega,s,y_s,z_s)ds-\int_t^{\widetilde\tau} z_sd\widetilde W_s,\quad 0\le t\le\widetilde\tau.
\enq
Conversely, if $(y_t,z_t)$ is a solution of \eqref{eq:eq36}, $(Y_t,Z_t):=(y_{\phi(t)},z_{\phi(t)}\cdot (\phi '(t))^{1/2})$ is a solution of \eqref{eq:eq31} using that $(\phi^{-1})'(\phi(t))=(\phi')^{-1}(t)$.
As in Section \ref{sec2}, $\widetilde f$ satisfies the uniform Lipschtz continuity. In fact, noting that $(\phi^{-1})'(t)=\alpha^{-2}(\phi^{-1}(t))$,
\begin{align}\label{eq:eq37}
\|\widetilde f(\omega,t,y,z)-&\widetilde f(\omega,t,y',z')\| \le\alpha^{-2}(\phi_t^{-1})[r_{\phi^{-1}(t)}\|y-y'\|+u_{\phi^{-1}(t)}\cdot\alpha(\phi_t^{-1})\|z-z'\|]\nonumber\\
&=\frac{r_{\phi^{-1}(t)}}{\max\{r_{\phi^{-1}(t)},u_{\phi^{-1}(t)}^2\}}\|y-y'\|+\frac{u_{\phi^{-1}(t)}}{\sqrt{\max\{r_{\phi^{-1}(t)},u_{\phi^{-1}(t)}^2\}}}\|z_t-z'_t\|\nonumber\\
&\le\|y-y'\|+\|z-z'\|,\quad dt\times d\P-a.s.
\end{align}
Now we are prepared to state some results on Wiener type BSDEs.
\begin{lem}\label{lem:lem31} Let the following conditions hold for BSDE \eqref{eq:eq31}.
\begin{enumerate}
\item The stochastic Lipschtz condition \eqref{eq:eq32} holds.
\item $\xi\in L_{\rho}^2 (\tau;\phi)$ and $f(\omega,s,0,0)/\alpha(s)\in L_{\rho}^2(0,\tau,\phi)$ for some $\rho>3$.  
\end{enumerate}
Then BSDE \eqref{eq:eq31} has a unique solution $(Y_t,Z_t)$ in $M_3^{2,\alpha}(0,\tau;\phi;\R^k;\R^{k\times d})$.
This solution actually belongs to $M_{\rho}^{2,\alpha}(0,\tau;\phi;\R^k;\R^{k\times d})$ and $Y\in U_{\rho}^2(0,\tau;\phi)$.
\end{lem}
\begin{proof}
We can notice that $\xi\in \widetilde L_{\rho}^2 (\widetilde\tau)$ and $\widetilde f(\omega,s,0,0)\in \widetilde L_{\rho}^2(0,\widetilde\tau)$ from
\begin{align}
&\int_0^{\tau}\exp(\rho\phi (s))\bigg\|\frac{f(\omega,s,0,0)}{\alpha(s)}\bigg\|^2ds=\int_0^{\phi (\tau)}\exp(\rho s)\biggl\|\frac{f(\omega,\phi^{-1}(s),0,0)}{\alpha (\phi^{-1}(s))}\biggr\|^2d\phi^{-1}(s)\nonumber\\
&=\int_0^{\phi (\tau)}\exp(\rho s)\|f(\omega,\phi^{-1}(s),0,0)\|^2 \bigl [(\phi^{-1})'(s)\bigr ]^2 ds=\int_0^{\widetilde\tau}\exp(\rho s)\|\widetilde f (\omega,s,0,0)\|^2 ds\nonumber.
\end{align}
Recalling \eqref{eq:eq37}, the simple application of \cite{DP}, Theorem 3.4 admits that BSDE \eqref{eq:eq36} has a unique solution $(y_t,z_t)$ in $\widetilde L_{3}^2(0,\widetilde\tau;\R^k\times\R^{k\times d})$ which belongs to $\widetilde L_{\rho}^2(0,\widetilde\tau;\R^k\times\R^{k\times d})$ and $y\in\widetilde U_{\rho}^2(0,\widetilde\tau)$.
Noting that $(Y_t,Z_t)=(y_{\phi(t)},z_{\phi(t)}\cdot[(\phi'(t)]^{1/2})$, we get the following expressions:
\begin{align}
\E\biggl[\int_0^{\widetilde\tau}\exp(\rho s)\|y_s\|^2ds\biggr]&=\E\biggl[\int_0^{\tau}\exp(\rho\phi(s))\|y_{\phi(s)}\|^2\phi'(s)ds\biggr]\nonumber\\
&=\E\biggl[\int_0^{\tau}\exp(\rho\phi(s))\|(\alpha\cdot Y)_s\|^2ds\biggr],\nonumber\\
\E\biggl[\int_0^{\widetilde\tau}\exp(\rho s)\|z_s\|^2ds\biggr]&=\E\biggl[\int_0^{\tau}\exp(\rho\phi(s))\big\|z_{\phi(s)}[\phi'(s)]^{1/2}\big\|^2ds\biggr]\nonumber\\
&=\E\biggl[\int_0^{\tau}\exp(\rho\phi(s))\|Z_s\|^2ds\biggr],\nonumber
\end{align}
\[
\E\big[\sup\{\exp{(\rho s)}\|y_s\|^2:0\le s\le \widetilde\tau\}\big]=\E\big[\sup\{\exp{\rho \phi(s)}\|Y_s\|^2:0\le s\le \tau\}\big].
\]
These are sufficient to complete the proof.
\end{proof}
In the above lemma, the stochastic Lipschtz condition in $y$ can be relaxed whereas BSDEs with random terminal time are well adopted under the monotonicity condition. This naturally admits us to give the following main result.
\begin{thm}\label{thm:thm32}
Suppose that the following conditions hold for BSDE \eqref{eq:eq31}.
\begin{enumerate}
\item There exist non-negative progressive processes $u_t,l_t$ and progressive process $r_t$ satisfying
\par $\alpha^2(t):=\max\{r_t^-,l_t,u_t^2\}>\epsilon$ ($r^-:=\max\{-r,0\}$) for some $\epsilon>0$ such that for any $y,y'\in\R^k$ and $z,z'\in\R^{k\times d}$;
\begin{description}
\item{1.1.} $\|f(\omega,t,y,z)\|\le \|f(\omega,t,0,z)\|+l_t(\|y\|+l')$ (where $l'\in{\{0,1\}}$)
\item{1.2.} $(y-y')(f(\omega,t,y,z)-f(\omega,t,y',z'))\le -r_t\|y-y'\|^2$
\item{1.3.} $\|f(\omega,t,y,z)-f(\omega,t,y',z')\|\le u_t\|z-z'\|$
\end{description}
\item $f(\omega,t,\cdot,z)$ is continuous.
\item $(\xi+l')\in L_{\rho}^2 (\tau;\phi)$ and $f(\omega,s,0,0)/\alpha(s)\in L_{\rho}^2(0,\tau;\phi)$ for some $\rho>3$,
 where $\phi(t):=\int_0^t\alpha^2(s)ds$.
\end{enumerate}
Then the conclusion of Lemma \ref{lem:lem31} holds.
\end{thm}
\begin{proof}
It can be easily seen that $\widetilde f$ is Lipschtz continuous in $z$.  The monotonicity and linear growth in $y$ are shown as follows.  
\begin{align}
(y-y')(\widetilde f(\omega,t,y,z)-&\widetilde f(\omega,t,y',z)) \le\alpha^{-2}(\phi_t^{-1})\cdot (-r_{\phi^{-1}(t)})\|y-y'\|^2\nonumber\\
&\le\alpha^{-2}(\phi_t^{-1})\cdot r_{\phi^{-1}(t)}^- \|y-y'\|^2=\frac{r_{\phi^{-1}(t)}^-}{\max\{r_{\phi^{-1}(t)}^-,l_{\phi^{-1}(t)},u_{\phi^{-1}(t)}^2\}}\|y-y'\|^2\le\|y-y'\|^2\nonumber.
\end{align}
\begin{align}
\|\widetilde f(\omega,t,y,z)\|&\le \|\widetilde f(\omega,t,0,z)\|+\alpha^{-2}(\phi_t^{-1})\cdot l_{\phi^{-1}(t)}(\|y\|+l')\nonumber\\
&=\frac{l_{\phi^{-1}(t)}}{\max\{r_{\phi^{-1}(t)}^-,l_{\phi^{-1}(t)},u_{\phi^{-1}(t)}^2\}}(\|y\|+l')\le \|\widetilde f(\omega,t,0,z)\|+(\|y\|+l')\nonumber
\end{align}
Now, the result easily follows from \cite{DP}, Theorem 3.4. 
\end{proof}

\begin{rem}\label{rem:rem31}
Existence and uniqueness results for BSDEs driven by Brownian motion with stochastic Lipschtz coefficients or stochastic monotone coefficients were already given in \cite{BEN, BK} under stronger assumptions than ours on linear growth coefficient and $\rho$ need to be enough large. For example, $\rho$ is assumed to be larger than 90 in \cite{ElH}(see the proof of Theorem 6.1 therein).
\end{rem}
For the Wiener-type BSDE with random terminal time, the comparison theorem also holds under the conditions for the existence and uniqueness (see \cite{DP}, Corollary 4.4.2).
Thus it is trivial that the comparison theorem holds for the BSDE \eqref{eq:eq31}.
Here we give the stability with respect to perturbations. Comparing to Theorem 3 in \cite{BK} we study under weaker assumptions.
\begin{thm}\label{thm:thm33}
Suppose $(\tau,\xi,f),(\tau',\xi',f')$ are the triples verifying the assumptions of    Theorem \ref{thm:thm32} with the same $\rho>3$. Let  $\Delta Y:=Y-Y',\Delta Z:=Z-Z'$ for $(Y,Z)\in L_{\rho}^{2,\alpha}(0,\tau;\phi;\R^k\times\R^{k\times d})$ and $(Y',Z')\in L_{\rho}^{2,\alpha}(0,\tau';\phi;\R^k\times\R^{k\times d})$ which are solutions of \eqref{eq:eq31} corresponding to $(\tau,\xi,f),(\tau',\xi',f')$, respectively.
Then there exist positive numbers $\beta,\delta$ for $3<\theta\le \rho$ such that
\begin{align}
\|\Delta Y(0)&\|^2+\beta\E\biggl[\int_0^{\tau\vee\tau'}\exp(\theta\phi(s))\alpha^2(s)\big(\|\Delta Y(s)\|^2+\|\Delta Z(s)\|^2\big)ds\biggr]\nonumber\\
&\le\E\|\exp(\theta\phi(\tau)/2)\xi-\exp(\theta\phi(\tau')/2)\xi'\|^2\nonumber\\
&\quad+\delta^{-1}\E\biggl[\int_0^{\tau\vee\tau'}\exp(\theta\phi(s))\biggl\|\frac{f(\omega,s,Y(s),Z(s))-f'(\omega,s,Y(s),Z(s))}{\alpha(s)}\biggr\|^2ds\biggr]\nonumber.
\end{align}
\end{thm}
\begin{proof}
We can adopt the same strategy as the proof of Theorem \ref{lem:lem31} thanks to \cite{DP}, Theorem 4.4 and omit the proof. 
\end{proof}
We can consider the case where the driver satisfies stochastic polynomial condition, that is, conditon 1.1 in Theorem \ref{thm:thm32} can be replaced by
\[
\|f(\omega,t,y,z)\|\le\|f(\omega,t,0,z)\|+l_t(\|y\|^p+l'), \quad p>1
\] 
or more generally
\beq\label{eq:eq38}
\|f(\omega,t,y,z)\|\le\|f(\omega,t,0,z)\|+l_t(\psi(\|y\|)+l')
\enq
for some continuous, increasing function $\psi$.
In this case, we can refer to \cite{BCa} (or \cite{Pa}) and there will not be any difficulty.
On the other hand, if the stochastic monotone coefficient is always non-negative (that is strictly monotone), we can set $r_t=0$, so it is sufficient to suppose that $\rho>1$ in preceding results.
If the driver is monotone decreasing, we can get a more useful result by referring to \cite{R}.
\begin{thm}\label{thm:thm34}
For BSDE \eqref{eq:eq31}, we suppose that the conditions 1.2, 1.3 and 2 in    Theorem \ref{thm:thm32} and \eqref{eq:eq38} hold with $l_t=l,r_t=0,k=1,l'=0$ for some $l,r\in\R$. 
\par We further assume that $\forall t\ge 0,f(\omega,t,0,0)=0$ and $|\xi|\le M$ for some $M\ge 0$.
Then there exists a solution $(Y_t,Z_t)$ of \eqref{eq:eq31} such that $|Y|\le M$ and $\forall t\ge 0,\E\biggl[\int_0^{t\wedge\tau}\|Z_s\|^2\biggr]<\infty$. 
\end{thm}
\begin{proof}
We only sketch the proof. We define a process $\phi(t):=\int_0^t (u^2(s)+1)ds$ with which we associate time change. Obviously $\phi^{-1}(t)\le t$. Let $\widetilde f$ denote the driver of time changed BSDE. 
\par Then it is easy to see that $\widetilde f$ is uniformly Lipschtz in $z$ with Lipsctz coefficient 1 and monotone decreasing in $y$. It also satisfies controlled growth condition with coefficient 1. We can easily check that $\widetilde f(\omega,t,0,0,)=0$.
\par So there exists a solution $(y_t,z_t)$ to the time changed BSDE such that $|y|<M$ and for any $t$, $\int_0^{t\wedge\widetilde\tau}\|z_s\|^2ds<\infty$ from \cite{R}, Theorem 3.1. Noting that $Y_t=y_{\phi(t)}$ and $\int_0^{t\wedge\tau}\|Z_s\|^2ds\le\int_0^{\phi^{-1}(t)\wedge\tau}\|Z_s\|^2ds=\int_0^{\phi^{-1}(t)\wedge\tau}\|Z_{\phi(s)}\|^2\phi'(s)ds=\int_0^{t\wedge\widetilde\tau}\|z_s\|^2ds<\infty$, we can complete the proof.
\end{proof}
\begin{rem}\label{rem:rem32}
We note that the uniqueness and comparison can be also stated under the further conditions using Theorems 3.6 and 3.7 in \cite{R}.
In    Theorem \ref{thm:thm34}, the exponential integrability condition on terminal value and the driver are not made and the same conditions as the case of uniformly Lipschtz were used for the study of the BSDE with stochastic one.
This is because the monotone coefficient which makes discounting rate is equal to zero.
\end{rem}


\subsection{Some aspects of further applications}\label{subsec3.1}
Concluding this section, we shall briefly mention that it is possible to have some further applications to get better or new results.
For example, it is not difficult to study $L^p-$solution \cite{BDHP}, the stability \cite{BH, T} and reflected BSDE \cite{AEN} in the context of  stochastic monotonicity condition through our framework. Although  $L^p-$solution of BSDE with stochastic Lipschtz condition was already studied in \cite{JQQ}, we can make the improved version, due to the preceding results.  We can use the results in \cite{Par} to study the BSDEs with jumps whose drivers are stochastic Lipschtz. 
For the martingale-type BSDE with stochastic Lipschtz condition, results of \cite{RW} are available.
The stochastic partial differential equations(SPDEs) with stochastic Lipschtz terms are connected to the backward doubly SDEs(BDSDEs) with stochastic Lipschtz coefficients and we can  refer to \cite{MS} concerning  the BDSDEs with random terminal time. Perhaps the derived results will be better than the ones in \cite{O1, O2} where the constant parameter appeared in integrability condition need to be sufficiently large.
The proposed technique can be also applied to the BSDEs on manifolds with “geometrical” Lipschtz condition studied in \cite{Bl}, which includes the results with respect to the random terminal time (see Section 5 therein), so the “geometrical” Lipschtz condition can be relaxed from the uniform one to the stochastic one.
The wellposedness of Mean-field backward stochastic delay equation with stopping terminal time and Lipschtz driver was stated in \cite{ML}, Theorem 3.1, so we can state the counterpart when the Lipschtz continuity is stochastic one. 
By referring to \cite{LRTY} where the results of second-order BSDEs (2BSDEs) with random terminal time are established, we can study 2BSDEs with stochastic Lipschtz condition.
The details are left to the readers and some of them may be non-trivial.
~\\
~\\
In this section, we used the results obtained under the monotonicity assumption.
So, the stronger integrability conditions on the driver and the solutions were still required as in the previous works.
We shall apply the proposed technique to the undiscounted BSDEs driven by Markov chains without the monotonicity assumption on driver in the next section.
\section{Markov chain BSDEs with stochastic Lipschtz coefficients}\label{sec4}
The BSDEs on Markov chains were first introduced in \cite{C2} and have developed in several papers, for example, the comparison theorem in \cite{C3} or the case of random terminal time in \cite{C1}. 
We present some preliminaries of the Markov chain BSDEs below.
\par Consider a continuous time, countable state Markov chain $X$ on $(\Omega,\F,\P,\mathbb F)$, where $\mathbb F$ is the natural filtration generated by $X$.
Without loss of generality, we assume that $X$ takes values from the unit vector $e_i$ in $\R^N,(N\in\N\cup\{\infty\})$, where $N$ is the number of states of the chain. We denote by $\Pi$ the state space.
If $A_t$ denotes the rate matrix of the chain at time $t$, then $(A_t)_{ij}\ge 0,\ i\neq j$ and $\forall j, \Sigma_i(A_t)_{ij}=0$.
For the simplicity, we shall assume that $A$ is uniformly bounded. The Markov chain $X$ has the following Doob-Meyer decomposition (see \cite{E}, Appendix B).
\beq\label{eq:eq41}
X_t=X_0+\int_0^t A_u X_{u-}du+M_t,
\enq
where $M$ is a pure discontinuous martingale with finite variation.
In this section, we further assume the Markov chain has the strong Markov property.
Let us consider the following BSDE to stopping time on Markov chain.
\beq\label{eq:eq42}
Y_t=\xi+\int_t^{\tau} f(\omega,u,Y_{u-},Z_u)du-\int_t^{\tau}Z_udM_u,\quad 0\le t\le\tau,
\enq
where $f:\Omega\times\R^+\times\R\times\R^N\rightarrow\R$ and $\int_0^t Z_udM_u=\sum_{i=1}^{N}\int_0^tZ_u^idM_u^i$.
Note that \eqref{eq:eq42} is contained in the class of BSDEs defined by \eqref{eq:eq14} due to \cite{C2}, Lemma 3.1 where it was shown that the sequence $(M^i),i=1,2,...$ has the martingale representation.
\begin{definition}\label{def:def41}
We define $\psi_t:=diag(A_tX_{t-})-A_tdiag(X_{t-})-diag(X_{t-})A_t^{\mathsf T}$.
Then the matrix $\psi_t$ is symmetric and positive (semi-)definite and $d<M>_t=\psi_tdt$ (see \cite{C2}).
\par Due to \eqref{eq:eq16}, we can set the stochastic semi-norm $\|\cdot\|_{M_t}$ as follows. 
\beq\label{eq:eq43}
\|Z\|_{M_t}^2:=Z_u^{\mathsf T}\psi_uZ_u,\quad Z\in \R^{1\times N}.
\enq
\end{definition}
We give further definitions from \cite{C1}.
\begin{definition}\label{def:def42}
We say that the driver $f$ is $\gamma-$balanced if there exists a random field $\eta:\Omega\times\R^+\times\R^N\times\R^N\rightarrow\R^N$, with $\eta(\cdot,\cdot,z,z')$ predictable and $\eta(\omega,t,0,0)$ Borel-measurable, such that 
\begin{itemize}
\item $f(\omega,t,y,z)-f(\omega,t,y,z')=(z-z')^{\mathsf T}(\eta(\omega,t,z,z')-AX_{t-})$
\item for each $e^i\in\Pi$, $(e_i^{\mathsf T}\eta(\omega,t,z,z'))/(e_i^{\mathsf T}AX_{t-})\in[\gamma,\gamma^{-1}]\ $ for some $\gamma>0$, where $0/0:=1$
\item $\textbf 1^{\mathsf T}\eta(\omega,t,z,z')=0\ $ for $\textbf 1\in\R^N$  the vector with all entries 1
\item $\eta(\omega,t,z+\alpha\textbf 1,z')=\eta(\omega,t,z,z')\ $ for all $\alpha\in\R$
\end{itemize}
\end{definition}
\begin{rem}\label{rem:rem41}
Note that the above definition of $\gamma-$balanced driver is closely connected to the notion of balanced driver in Section \ref{sec2} (see \cite{C1}, Lemma 3).
\end{rem}
\begin{definition}\label{def:def43}
Let $\mathcal Q_{\gamma}$ denote the family of all measures $Q$ where $X$ has the compensator $\eta(t,\omega)$, for $\eta$ a predictable process with $\textbf 1^{\mathsf T}\eta(t,\omega)=0$ and $\frac{e_i^{\mathsf T}\eta(t,\omega)}{e_i^{\mathsf T}AX_{t-}}\in[\gamma,\gamma^{-1}]$ for all $e_i\in\Pi$, where $0/0:=1$.
That is, $X_t=X_0+\int_0^t\eta_tdt+Q$-martingale,$\ Q\in\mathcal Q^{\gamma}$.
\end{definition}
We give the key result of \cite{C1} (see Theorem 3, Remark 4 therein).
\begin{lem}\label{lem:lem41}
Suppose that the following conditions are verified for Markov chain BSDE \eqref{eq:eq42}.
\begin{enumerate}
\item $\xi$ is $\F_{\tau}-$measurable.
\item There exist non-decreasing functions $K_1,K_2:\R^+\rightarrow[1,\infty)$ and some constants $\beta,\widetilde\beta>0$ such that
\[
\E^Q[\xi|\F_t]\le K_1(t),\E^Q[(1+\tau)^{1+\beta}|\F_t]\le K_1(t), \E^Q[K_1(\tau)^{1+\widetilde\beta}|\F_t]\le K_2(t),
\] 
for all $\P-a.s.$ all $Q\in\mathcal Q_{\gamma}$ and all $t$.
\item $f:\Omega\times\R^+\times\R\times\R^N\rightarrow\R$ is $\gamma-$balanced.
\item The discounting terms are uniformly bounded above, that is, there exists a constant $C_1\in\R$ such that for any $y,y',z$ and $s<t$,
\[
\int_s^t r(\omega,u,y,y',z)du<C_1,\quad r(\omega,u,y,y',z):=\frac{f(\omega,t,y,z)-f(\omega,t,y',z)}{y-y'}.
\]
\item There exists $C_2\in\R$, $\hat\beta\in[0,\beta]$ such that $|f(\omega,t,0,0)|\le C_2(1+t^{\hat\beta})$.
\item $f:\Omega\times\R^+\times\R\times\R^N\rightarrow\R$ is uniformly Lipschtz in $y$. That is, there exists a constant $C$ such that $|f(\omega,t,y,z)-f(\omega,t,y',z)|\le C|y-y'|$ for all $y,y',z$.
\end{enumerate}
Then the BSDE \eqref{eq:eq42} has a unique solution such that $|Y_t|\le (1+C_2)\exp(C_1)|K_1(t)|$.
\end{lem}
\begin{rem}\label{rem:rem42}
Note that the fourth condition is verified if the driver is monotone decreasing.
\end{rem}
In  Lemma \ref{lem:lem41}, the Lipschtz condition in $y$ is restrictive.
We give a simple illustration below with the motion of a particle on graph.
Consider a model for transmission of messages from a node to another node over a network. Let the chain $X$ describe the motion of a message. Then the probability that the message reaches its target is given as the solution of the following BSDE (see \cite{C1}, Section 4).
\beq\label{eq:eq44}
I_{\{X_{\tau}=x_1\}}=Y_t-\int_t^\tau-r_{X_{u-}}Y_{u-}du+\int_t^\tau Z_u dM_u,\quad 0\le t\le\tau,
\enq
where $r_x$ is the rate by which the node $x$ loses a message. To suppose that the losing rate at each node is bounded is an assumption rarely satisfied in real world. It depends on the time variable in general and it should be written as $r(t,X_{t-})$ which may be unbounded.

\par The main result of this section is as follows (we shall give the proof later).
\begin{thm}\label{thm:thm42}
Suppose that the conditions 1)$-$5) in  Lemma \ref{lem:lem41} are satisfied for BSDE \eqref{eq:eq42}. Let the driver $f$ be stochastic Lipschtz in $y$, that is, there exists a non-negative predictable process $C(t)$ such that for any $t,y,z,z'$,
\beq\label{eq:eq45}
|f(\omega,t,y,z)-f(\omega,t,y',z)|\le C(t)|y-y'|.
\enq
Then BSDE \eqref{eq:eq42} has a unique solution such that $|Y_t|\le 2\exp(C_1)|K_1(t)|$.
\end{thm}
\begin{rem}\label{rem:rem43}
Theorem \ref{thm:thm42} admits that the bound of a solution does not depend on $C_2$. From this fact, the bound estimate on a solution  can be replaced by $|Y_t|\le(1+C_2\wedge1)\exp(C_1)|K_1(t)|$ in  Lemma \ref{lem:lem41}.
\end{rem}
Usually, if one wants to relax the Lipschtz continuity as the stochastic one, it has to be considered that the stronger integrability conditions on terminal value, driver and the solution are required, instead. However, this is not true for undiscounted BSDE.
It is because the terminal value and driver of this BSDEs are not needed to be discounted at some rate and one can consider the direct conditions on them respectively.
~\\
\par In Theorem \ref{thm:thm42}, the conditions on stopping time seem to be unfamiliar and it is required to afford an example when they are satisfied. In this context, S. N. Cohen \cite{C1} showed that the direct conditions on stopping time are satisfied when the stopping time is a hitting time of a subset of $\Pi$ under the uniform ergodicity of the chain by the way of examining the exponential ergodicity of the chain under the perturbations of rate matrix. One can observe that the above hitting time only depends on the character of the chain. In    Theorem \ref{thm:thm42}, the driver is stochastic Lipschtz only in $y$ and the $\gamma-$balanced condition related to $z$ is still required.
On the other hand, it was shown  in two uniform and stochastic Lipschtz settings that the conditions on stopping time and terminal value for the wellposedness of BSDE \eqref{eq:eq42} coincide. These lead to the following result (see \cite{C1}, Lemma 6).

\begin{lem}\label{lem:lem44}
Suppose that rate matrix is time-homogeneous under the measure $\P$ and the chain is uniformly ergodic.
Let $\tau$ be the first hitting time of a set $\Xi\subseteq\Pi$ and $\xi$ be a random variable of the form $\xi=g(\tau,X_{\tau})$ for some function $g(t,x)\le k(1+t^{\beta})$ for some $k,\beta>0$. Then there exist functions $K_1,K_2$ satisfying the requirements of    Theorem \ref{thm:thm42}.
\end{lem}

When the terminal time and terminal value have the forms like in  Lemma \ref{lem:lem44} and the driver is Markovian, that is, $f(\omega,t,y,z)=\bar f(X_{t-},t,y,z)$ for some $\bar f$, we can give the ODE system with boundary condition which describes the solution of BSDE in the context of stochastic Lipschtz assumption (see \cite{C1}, Theorems 6 or 7).   
Now we seem to prove Theorem \ref{thm:thm42} by means of time change described in Section \ref{sec2}. 
~\\

\begin{profof}\textbf{Thoerem 4.2}
We define the process $\phi$ as follows (this is based on the same idea as in Section \ref{sec2}).
\beq\label{eq:eq46}
\phi(t):=\int_0^t \alpha^2(s)ds,\quad \alpha^2(s):=\max\{C(s),C_2,1\}.
\enq
Then it follows that $t\ge\phi^{-1}(t)$ from $\phi(t)\ge\int_0^t1ds=t$. We set $\widetilde\F_t:=\F_{\phi^{-1}(t)},\widetilde{\mathbb F}:=\{\widetilde{\F}_t\}_{t\ge 0}$ as in Section \ref{sec2}.
As $X$ is a strong Markov chain, $\widetilde X:=X_{\phi^{-1}(t)}$ is also a strong Markov chain with respect to $\widetilde{\mathbb F}$ (e.g. see \cite{Ba}, Chapter 22, Section 3). 
~\\
Using the expression \eqref{eq:eq41},
\beq\label{eq:eq47}
\widetilde X_t=X_0+\int_0^{\phi^{-1}(t)}A_uX_udu+M_{\phi^{-1}(t)}=
\widetilde X_0+\int_0^t \widetilde A_u \widetilde X_u du+\widetilde M_t,
\enq
where $\widetilde A_u:=A_{\phi^{-1}(u)}\cdot(\phi^{-1})'(u)$ and $\widetilde M_t:=M_{\phi^{-1}}(t)$.
We recall that $\widetilde M=(\widetilde M^i),i=1,2,...,N$ is a sequence of orthogonal martingales which has martingale representation on $(\Omega,\F,\P,\widetilde{\mathbb F})$ (see  Lemma \ref{lem:lem24}).
Therefore, \eqref{eq:eq47} is the (unique) Doob-Meyer decomposition of $\widetilde X$.
If we denote by $R_u$ the rate matrix of $\widetilde X$,  then 
\beq\label{eq:eq48}
R_u\widetilde X_u=\widetilde A_u\widetilde X_u,\quad dt\times d\P-a.s.\ .  
\enq
It follows that $R_u$ is uniformly bounded from $(\phi^{-1})'(t)=\frac{1}{\phi'(\phi_t^{-1})}=\alpha^{-2}(\phi^{-1}(t))\le 1$, so the $\widetilde{\mathcal F}-$chain $\widetilde X$ is also regular.
We can consider that the random rate matrix $\widetilde A_t$ plays the role of transiation rate matrix of $\widetilde X$.
Next, we shall show that $\widetilde f(\omega,s,y,z):=f(\omega,\phi^{-1}(s),y,z)\cdot (\phi^{-1})'(s)$ is $\gamma-$balanced with respect to $\widetilde{\mathbb F}$. We define $\widetilde\eta(\omega,t,z,z'):=\eta(\omega,\phi^{-1}(t),z,z')\cdot (\phi^{-1})'(t)$. Then by the definition and \eqref{eq:eq48}, we have
\begin{align}
\widetilde f(\omega,t,y,z)-\widetilde f(\omega,t,y,z')&=(f(\omega,\phi^{-1}(t),y,z)-f(\omega,\phi^{-1}(t),y,z')) (\phi^{-1})'(t)\nonumber\\
&=(z-z')^{\mathsf T}(\eta(\omega,\phi^{-1}(t),z,z')-A_{\phi^{-1}(t)}X_{\phi^{-1}(t)}) (\phi^{-1})'(t)\nonumber\\
&=(z-z')^{\mathsf T}(\widetilde\eta(\omega,t,z,z')-\widetilde A_t\widetilde X_t)=(z-z')^{\mathsf T}(\widetilde\eta(\omega,t,z,z')-R_t\widetilde X_t),\nonumber\\ 
(e_i^{\mathsf T}\widetilde{\eta}(\omega,t,z,z'))/(e_i^{\mathsf T}R_t\widetilde X_t)&=(e_i^{\mathsf T}\widetilde{\eta}(\omega,t,z,z'))/(e_i^{\mathsf T}\widetilde A_t\widetilde X_t)\nonumber\\
&=(e_i^{\mathsf T}\eta(\omega,\phi^{-1}(t),z,z'))/(e_i^{\mathsf T}A_{\phi^{-1}(t)} X_{\phi^{-1}(t)})\in[\gamma,\gamma^{-1}],\nonumber\\
\textbf{1}^{\mathsf T}\widetilde\eta(\omega,t,z,z')&=\textbf{1}^{\mathsf T}\eta(\omega,\phi^{-1}(t),z,z')(\phi^{-1})'(t)=0,\nonumber\\
\widetilde\eta(\omega,t,z+\alpha\textbf{1},z')&=\eta(\omega,\phi^{-1}(t),z+\alpha\textbf{1},z')\cdot (\phi^{-1})'(t)\nonumber\\
&=\eta(\omega,\phi^{-1}(t),z,z')\cdot (\phi^{-1})'(t)=\widetilde\eta(\omega,t,z,z'),\nonumber
\end{align}
$dt\times d\P-a.s.$
\par So $\widetilde f$ is $\gamma-$balanced.
We note that $\widetilde f$ is uniformly Lipschtz in $z$ under norm $\|\cdot\|_{\widetilde M_t}$ because it is $\gamma-$balanced (see \cite{C1}, Lemma 1).
From the expressions \eqref{eq:eq47} and \eqref{eq:eq48}, it is trivial that the family of probability measures where $\widetilde X$ has the $\widetilde{\mathbb F}-$predictable compensator $\widetilde\eta(t,\omega)$ such that $\textbf{1}^{\mathsf T}\widetilde\eta=0$ and $\forall 1\le i\le N; (e_i^{\mathsf T}\widetilde\eta(t,\omega))/(e_i^{\mathsf T}R_t\widetilde X_{t-})\in[\gamma,\gamma^{-1}]$ is also $\mathcal Q_{\gamma}$.  
Finally we show that the following time changed BSDE has a unique solution($\widetilde\tau:=\phi(\tau)$).
\beq\label{eq:eq49}
y_t=\xi+\int_t^{\widetilde\tau}\widetilde f(\omega,s,y_s,z_s)ds-\int_t^{\widetilde\tau}z_sd\widetilde M_s.
\enq
We have already seen that $\widetilde f$ is $\gamma-$balanced.
\par Let us define the non-decreasing functions $\widetilde K_1(t):=K_1(\phi^{-1}(t))$ and $\widetilde K_2(t):=K_2(\phi^{-1}(t))$. Then $\forall Q\in\mathcal Q_{\gamma};\E^Q[\xi|\widetilde{\F}_t]\le K_1(\phi^{-1}(t))=\widetilde K_1(t)$ and $\E^Q[\widetilde K_1(\widetilde{\tau})^{1+\widetilde{\beta}}|\widetilde{\F}_t] =\E^Q[K_1(\tau)^{1+\widetilde{\beta}}|\F_{\phi^{-1}(t)}]\le K_2(\phi^{-1}(t))=\widetilde K_2(t)$. 
\par Using the assumptions on $f$, we can get the following expressions on $\widetilde f$. 
\begin{align}
|\widetilde f(\omega,t,0,0)|&=\alpha^{-2}(\phi^{-1}(t))|f(\omega,\phi^{-1}(t),0,0)|\nonumber\\
&\le(1+\phi^{-1}(t)^{\hat{\beta}})\cdot C_2/[C(\phi^{-1}(t))+C_2+1]\le 1+\phi^{-1}(t)^{\hat{\beta}} \le1+t^{\hat{\beta}}\nonumber,
\end{align} 
\begin{align}
\int_s^t (\widetilde f(\omega,u,y,z)-\widetilde f(\omega,u,y',z))/(y-y')du&=\int_s^t r(\omega,\phi^{-1}(u),y,y',z)d\phi^{-1}(u)\nonumber\\
&=\int_{\phi^{-1}(s)}^{\phi^{-1}(t)} r(\omega,u,y,y',z)du\le C_1\nonumber,
\end{align}
\begin{align}
|\widetilde f(\omega,t,y,z)-\widetilde f(\omega,t,y',z)|&=|\widetilde f(\omega,\phi^{-1}(t),y,z)-\widetilde f(\omega,\phi^{-1}(t),y',z)|\cdot\alpha^{-2}(\phi^{-1}(t))\nonumber\\
&\le C(\phi^{-1}(t))|y-y'|/\big(\max\{C(\phi^{-1}(t)),C_2,1\}\big)\le|y-y'|\nonumber.
\end{align}
So BSDE \eqref{eq:eq49} has a unique solution satisfying $|y_t|\le 2\exp(C_1)|\widetilde K_1(t)|$ by  Lemma \ref{lem:lem41}.
If we set $(Y_t,Z_t):=(y_{\phi(t)},z_{\phi(t)})$,    Theorem \ref{thm:thm25} shows that it is a solution of BSDE \eqref{eq:eq42}. Because the solution $y_t$ of \eqref{eq:eq49} is unique up to indistinguishability, $Y_t$ is also unique up to indistinguishability. And $|Y_t|=|y_{\phi(t)}|\le2\exp(C_1)|\widetilde K_1(\phi(t))|=2\exp(C_1)|K_1(t)|$.
\end{profof}
\begin{rem}\label{rem:rem45}
We note that the comparison theorem for BSDE \eqref{eq:eq42} holds under the stochastic Lipschtz condition from the corresponding comparison theorem for BSDE \eqref{eq:eq49} (see \cite{C1}, Theorem 5). 
\end{rem}

\subsection{Additional use of time change}\label{subsection41}
In Theorem \ref{thm:thm42}, the $5^{th}$ condition is not strictly necessary. In fact, it is sufficient to suppose that $f(\omega,s,0,0)/(1+t^{\hat\beta})$ is integrable on every finite interval in $\R^+$. 
For any $m>1$, consider the process $\bar{\phi}(t):=m\int_0^t(|f(\omega,s,0,0)|/(1+t^{\hat\beta})+1)ds$ with which we associate time change. Then
\[
|\widetilde f(\omega,t,0,0)|=\frac{|f(\omega,\bar\phi^{-1}(t),0,0)|}{1+|f(\omega,\bar\phi^{-1}(t),0,0)|/(1+\bar\phi^{-1}(t)^{\hat\beta})}\cdot\frac{1}{m}\le\frac{1+\bar\phi^{-1}(t)^{\hat\beta}}{m}\le\frac{1+t^{\hat\beta}}{m}
\]
Also we have the bound on solution such that $|Y_t|\le(1+1/m)\exp(C_1)|K_1(t)|,$ for all $m>1$. Taking $m\rightarrow\infty$, we get $|Y_t|\le\exp(C_1)|K_1(t)|$. So we can show the improvement of  Theorem \ref{thm:thm42}.
\begin{thm}\label{thm:thm44}
Suppose that conditions 1$-$4,\ 6 in Lemma \ref{lem:lem41} and the stochastic Lipschtz condition \eqref{eq:eq45} hold. We further assume that   $f(\omega,s,0,0)/(1+t^{\hat\beta})$ is integrable on every finite interval in $\R^+$.
Then there exists a unique solution to BSDE \eqref{eq:eq42} such that $|Y_t|\le\exp(C_1)|K_1(t)|$. If the driver is monotone decreasing, then $|Y_t|\le|K_1(t)|$.
\end{thm}
In this subsection, we made the use of time change away from the discussion on Lipschtz continuity. Perhaps, there will be other problems to which we can apply time change effectively in the range of stochastic calculus.
\section{Conclusion}\label{sec5}
In this paper, we showed that the technique for dealing with the BSDEs with stochastic Lipschtz coefficients by time change. The technique says that when we study the BSDE to stopping time, the Lipschtz condition can be given as the stochastic one.
Also roughly speaking, most of the results of BSDEs obtained under the Lipschtz continuity may be extended to the case of stochastic Lipschtz continuity.
Of course, this is only possible when we are aware of the results with respect to random terminal time and this admits the importance on the study of them.


\begin{thebibliography}{54}
\bibitem{AEN} A. Aman, A. Elouaflin and M. N’zi, Reflected generalized BSDEs with random time and applications, arXiv preprint, 2010, arXiv:1011.3223v1 [math.PR]
\bibitem{AIP} S. Ankirchner, P.Imkeller and A.popier, On measure solutions of backward stochastic  differential  equations, Stochastic Process. Appl. 119 (2009) 2744-2772. 
\bibitem{BEN} K. Bahlali, A. Elouaflin and M. N’zi, Backward stochastic differential equations with stochastic monotone coefficients, J. Appl. Math. Stoch. Anal. 4 (2004) 317-335. 
\bibitem{Ba} R. F. Bass, Stochastic Processes. Cambridge University Press, 2011.
\bibitem{BK} C. Bender and M. Kohlmann, BSDEs with Stochastic Lipschitz Condition, Technical report, Center of Finance and Econometrics, University of Konstanz, 2000. Available at http://hdl. handle.net/10419/85163.
\bibitem{Bi} J. M. Bismut, Conjugate  convex functions in optimal stochastic control, J. Math. Anal. Appl. 44 (1973) 384-404.
\bibitem{Bl} F. Blache,  Backward Stochastic Differential Equations on Manifolds, arXiv preprint, 2005, arXiv:math/0501265v1 [math.PR]
\bibitem{BC} P. Briand and F. Confortola, BSDEs with stochastic Lipschitz condition and quadratic PDEs in Hilbert spaces, Stochastic Process. Appl. 118(5) (2008) 818-838.
\bibitem{BH} P. Briand and Y. Hu, Stability of BSDEs with random terminal time and homogenization of semilinear elliptic PDEs, J. Funct. Anal.  155 (1998) 455-494.
\bibitem{BCa} P. Briand and R. Carmona, BSDEs with polynomial growth generators, J. Appl. Math. Stochastic Anal. 13 (2000) 207-238.
\bibitem{BDHP} P. Briand, B. Delyon, Y. Hu, E. Pardoux and L. Stoica, $L^p$ solutions of Backward  Stochastic Differential Equations, Stochastic Process. Appl. 108 (2003) 109-129.
\bibitem{BLM} P. Briand, J. P. Lepeltier, J. San Martín, One-dimensional BSDE’s whose coefficient is monotonic in $y$ and non-lipschitz in $z$, Bernoulli. 3(1) (2007) 80-91.
\bibitem{CFS} R. Carbone, B. Ferrario and M. Santacroce, Backward stochastic differential equations driven by c\`adl\`ag martingales, Theory Probab. Appl. 52(2) (2008) 304-314.
\bibitem{C1} S. N. Cohen, Undiscounted Markov chain BSDEs to stopping times, J. Appl. Probab. 51 (2014) 262-281.
\bibitem{C2} S. N. Cohen and R. J. Elliott, Solutions of backward stochastic differential equations  on Markov chains, Commun. Stoch. Anal.  2(2) (2008) 251-262.  
\bibitem{C3} S. N. Cohen and R. J. Elliott, Comparisons for backward stochastic  differential  equations on Markov chains and related no-arbitrage condition, Ann. Appl. Probab. 20(1) (2010)  267-311. 
\bibitem{C4} S. N. Cohen and R. J. Elliott, Existence, uniqueness and comparisons for BSDEs in general spaces, Ann. Probab. 40(5) (2012) 2264-2297.
\bibitem{C5} S. N. Cohen and R.J. Elliott, Stochastic Calculus and Applications, Springer New York, 2015.
\bibitem{DP} R. W. R. Darling and E. Pardoux, Backwards SDE with random terminal time and applications to semilinear elliptic PDE, Ann. Probab. 25 (1997) 1135-1159.
\bibitem{DV} M. H. A. Davis and P. Varaiya, On the multiplicity of an increasing family of sigma-fields, Ann. Probab. 2 (1974) 958-963.
\bibitem{D} M. H. A. Davis, Martingale Representation and All That, Advances in Control, Communication Networks and Transportation Systems: In Honor of Pravin Varaiya, E.H. Abed (Ed.), Systems and Control: Foundations and Applications Series, Birkhauser, Boston, (2005).
\bibitem{E} R. J. Elliott, L. Aggoun and J. B. Moore, Hidden Markov Models: Estimation and  Control. Springer-Verlag, Berlin-Heidelberg-New York, 1994.
\bibitem{ElH} N. El Karoui and S. J. Huang, A general result of existence and  uniqueness of  backward  stochastic  differential  equations,  In Backward Stochastic Differential   Equations, Pitman Research Notes in Mathematics Series. 364 (1997) 27-36. Longman, Harlow.
\bibitem{ElPQ} N. El Karoui, S. Peng and M.C. Quenez, Backward Stochastic Differential Equations in Finance, Mathematical Finance. 7 (1997) 1-71.
\bibitem{Hu} L. -Y. Hu, Reflected backward doubly stochastic differential equations driven by a L\`evy process with stochastic Lipschitz
condition, Appl. Math. Comput. 219 (2012) 1153–1157.
\bibitem{LM} J. P. Lepeltier and J. S. Martin, Backward stochastic differential equations with  continuous coefficients, Statist. Probab. Lett. 34 (1997) 425-430.
\bibitem{LRTY} Y. Lin, Z. Ren, N. Touzi and J. Yang, Second order backward SDE with random terminal time, arXiv preprint, 2018,  arXiv:1802.02260v1 [math.PR]
\bibitem{ML} H.P. Ma and B. Liu, Infinite horizon optimal control problem of mean-field backward stochastic delay differential equation under partial information, Eur. J. Control. 36 (2017) 43-50.
\bibitem{M} X. Mao, Adapted solution of backward stochastic differential equations with non-lipschitz coefficients, Stochastic Process. Appl.  58 (1995) 281-292.
\bibitem{MEl} M. Marzougue and M. El Otmani, Double barrier reflected BSDEs with stochastic Lipschitz  coefficient, arXiv preprint, 2018, arXiv:1801.01016v1 [math.PR].
\bibitem{MS} A. Matoussi, W. Sabbagh, Numerical Computation for backward doubly SDEs with random terminal time, arXiv preprint, 2016, arXiv:1409.2149v4 [math.PR]
\bibitem{Med} P. Medvegyev, Stochastic Integration Theory. Oxford University Press, 2006.
\bibitem{O1} J. M. Owo, Backward doubly stochastic differential equations with stochastic lipschitz condition, Statist. Probab. Lett. 96 (2015) 75-84.
\bibitem{O2} J. M. Owo, $L^p$-solutions of backward doubly stochastic differential equations with stochastic lipschtz condition and $p\in(1,2)$, ESAIM: PS. 21 (2017) 168-182.
\bibitem{PPS} A. Papapantoleon, D. Possamai and A. Saplaouras, Existence and uniqueness results for BSDEs with jumps: the whole nine yards, arXiv preprint, 2018, arXiv:1607.04214v3 [math.PR]
\bibitem{PP1} E. Pardoux and S. Peng, Adapted solution of a backward stochastic differential equation, Systems Control Lett. 1 (1990) 55-61.
\bibitem{PP2} E. Pardoux and S. Peng, Backward stochastic differential equations and quasilinear parabolic partial differential equations, in : B.L. Rozovskii, R.B. Sowers (eds) Stochastic partial differential equations and their applications (Lect. Notes Control Inf. Sci. 176, 200-217) Berlin, Heidelberg New York : Springer 1992.
\bibitem{Par}E.Pardoux, Generalized discontinuous backward stochastic differential equations, In Backward stochastic differential equations, Pitman Research Notes in Mathematics Series. 364 (1997) 207-219. Longman, Harlow.
\bibitem{Pa} E. Pardoux,  BSDEs, weak convergence and homogenization of semilinear PDEs, In: Nonlinear Analysis, Differential Equations and Control (ed. by F. H. Clarke and R. J. Stern), Kluwer Acad. Pub. (1999), 503-549.
\bibitem{PA} E. Pardoux, A. Rasscanu, Stochastic differential equations, Backward SDEs, Partial differential equations, Springer, Cham, (2014).
\bibitem{P2} S. Peng, Probabilistic interpretation for systems of quasilinear parabolic partial differential equations, Stochastics. 37 (1991) 61-74.
\bibitem{RY} D. Revuz and M. Yor, Continuous Martingales and Brownian  Motion. Springer, Berlin (1999).
\bibitem{R} M. Royer, BSDEs with a random terminal time driven by a monotone generator and their links with PDEs, Stochastics. 76 (2004) 281-307.
\bibitem{RW} F. Russo and L. Wurzer, Elliptic PDEs with distributional drift and backward SDEs driven by a c\`adl\`ag martingale with random terminal time, Stoch. Dyn. 17(4):1750030(36 pages) 2017.
\bibitem{T} S. Toldo, Stability of solutions of BSDEs with random terminal time, ESAIM: PS. 10 (2006) 141-163.
\bibitem{JQQ} J. Wang, Q. Ran and Q. Chen, $L^p$-Solutions of BSDEs with Stochastic Lipschitz Condition, J. Appl. Math. Stochastic Anal.
2007 (2006) 1-14.
\bibitem{W} L. Wen, Reflected BSDE with stochastic Lipschitz coefficient, arXiv preprint, 2009, arXiv:0912.2162v3 [math.PR]
\bibitem{YZ} J. M. Yong and X. Y. Zhou, Stochastic Controls: Hamiltonian Systems and HJB Equations, Springer, New York, 1999.
\end{thebibliography}
\end{document}